\documentclass[12pt]{article}
\usepackage{amssymb,amsmath,amsthm}
\usepackage{graphicx,epstopdf}
\usepackage{cite}
\usepackage{float}
\usepackage{tabularx}
\usepackage{graphicx}
\newtheorem{theorem}{Theorem}[section]

\newtheorem{lemma}{Lemma}[section]
\newtheorem{definition}{Definition}[section]
\newtheorem{remark}{Remark}[section]

\newcommand{\beq}{\begin{equation}}
	\newcommand{\eeq}{\end{equation}}
\newcommand{\beqn}{\begin{eqnarray}}
	\newcommand{\eeqn}{\end{eqnarray}}

\def\to{{\rightarrow}}

\allowdisplaybreaks \baselineskip 20pt

\newcommand{\dd}{\mathrm{d}}
\usepackage{color}
\begin{document}
	
	\title{Besicovitch-like Pseudo Almost Automorphic Stochastic Processes: Theory and Applications to Semilinear SDEs\thanks{This work is supported by the National Natural Science Foundation of China
under Grant No. 12261098.}}
	\author{Yongkun Li\thanks{The corresponding author. Email: yklie@ynu.edu.cn.} and Yuwei Cao \\
		Department of Mathematics, Yunnan University,\\
		Kunming, Yunnan 650091,  China}
	
	\date{}
	\maketitle{}
	\begin{abstract}
This paper introduces, for the first time, the concept of Besicovitch-like  pseudo almost automorphic stochastic processes. As a meaningful application, we investigate a class of semilinear stochastic equations driven by Brownian motion in a separable Hilbert space. By combining the Banach fixed point theorem with refined inequality techniques, we establish sufficient conditions guaranteeing the existence, uniqueness, and global exponential stability of Besicovitch-like  pseudo almost automorphic solutions in distribution for the considered equations. An example is provided to show the validity of the obtained results.
	\end{abstract}
	
	\allowdisplaybreaks
	\section{Introduction}
	\setcounter{equation}{0}
	\indent

Periodicity and almost periodicity are important objects of study in the qualitative theory of differential equations. Almost automorphy is a natural generalization of almost periodicity, and pseudo almost automorphy further generalizes almost automorphy. These concepts together constitute important qualitative behaviors for both deterministic and stochastic differential equations. Specifically, almost automorphic functions were first introduced by Bochner when he studied a problem in differential geometry \cite{bc}. Since then, this concept has been successively generalized in various senses, giving rise to notions of almost automorphy in senses such as Stepanov, Weyl, and Besicovitch \cite{2,b2,b3,b4,fyu}. It is worth noting that Weyl-like almost automorphy generalizes Stepanov-like almost automorphy, and Besicovitch-like almost automorphy further generalizes Weyl-like almost automorphy. In parallel, inspired by Zhang's concept of pseudo almost periodic functions \cite{zh1}, researchers have successively proposed concepts of pseudo almost automorphic functions in various senses, including pseudo almost automorphic functions, weighted pseudo almost automorphic functions, and various measure-based concepts of pseudo almost automorphic functions \cite{w1,w2,w3,w4,w5,w6,w7}.

Parallel to the deterministic framework, periodicity, almost periodicity, and almost automorphy of stochastic differential equations are also important topics in the qualitative study of stochastic differential equations. It should be pointed out here that, as noted in \cite{xx}, for stochastic differential equations it is more reasonable to study almost automorphic solutions in the sense of distribution than in the mean-square sense. Motivated by this, various notions of stochastic processes in the sense of distribution have been proposed and successfully applied to study the dynamical behavior of different types of stochastic differential equations. Recently, Ref. \cite{zl} introduced the concept of Stepanov-like almost automorphic stochastic processes in distribution. However, for the most complex class of almost automorphic oscillations, namely Besicovitch-like almost automorphy, the corresponding concept of Besicovitch-like pseudo almost automorphic stochastic processes has not yet been established.

Motivated by the above discussions, this paper introduces, for the first time, the concept of Besicovitch-like  pseudo almost automorphic stochastic processes. As an application of this new notion, we  study the existence, uniqueness and global exponential stability of  Besicovitch-like pseudo almost automorphic solutions for the following semilinear stochastic equation driven by Brownian motion in a separable Hilbert
space $\mathbb{X}$:
\begin{align}\label{1}
\mathrm{d}X(t)=AX(t)\mathrm{d}t+f(t,X(t))dt+g(t,X(t))\mathrm{d}W(t),\,\,t\in\mathbb R,
\end{align}
where $A :Dom(A )\subset \mathbb X \rightarrow  \mathbb{X}$ is the infinitesimal generator of a $C_0$-semigroup $T(t)$, functions $f:\mathbb R\times \mathcal{L}^p(\Omega,\mathbb{X}), \rightarrow \mathcal{L}^p(\Omega,\mathbb X)$ and $g: \mathbb{R}\times  \mathcal{L}^p(\Omega,\mathbb X)\rightarrow \mathbb{L}^0_2$, where $\mathbb{L}^0_2=\mathbb{L}_2^0(\mathbb{X},\mathbb{X})$ is a separable Hilbert space with respect to the Hilbert-Schmidt norm $\|\cdot\|_{\mathbb{L}^0_2}$  satisfying some additional conditions will be stated lader, and
$W(t)$ is a two-sided standard one-dimensional Brownian motion with values in $\mathbb{X}$ defined on the filtered probability space $(\Omega,\mathcal{F},\{\mathcal{F}_t\}_{t\in \mathbb{R}},P)$, where $\mathcal{F}_t$ satisfies the usual conditions.

The remainder of this paper is organized as follows. Section 2 introduces some necessary fundamental concepts and useful lemmas. Section 3 provides a rigorous definition of Besicovitch-like  pseudo almost automorphic stochastic processes. Section 4 establishes, by means of the fixed point method and inequality techniques, the existence, uniqueness, and global exponential stability of the desired solutions. Section 5 presents an example to verify the effectiveness of the theoretical results. Section 6 gives a brief conclusion.

	\section{Preliminaries}
	\setcounter{equation}{0}
	\indent

We assume that $(\mathbb{X},\|\cdot\|)$ and $(\mathbb{K},\|\cdot\|_\mathbb{K})$ are separable Hilbert spaces. Let $(\Omega, \mathcal{F}, \{\mathcal{F}_t\}, \mathbb{P})$ be a filtered probability space satisfying the usual conditions. The space $\mathbb{L}^0_2=\mathbb{L}_2(\mathbb{K},\mathbb{X})$ denotes the space of all Hilbert-Schmidt operators from $\mathbb{K}$ into $\mathbb{X}$, equipped with the Hilbert-Schmidt norm $\|\cdot\|_{\mathbb{L}^0_2}$; it is itself a Hilbert space.

 For $1\le p<\infty$, we denote by $\mathcal{L}^{p}(\Omega,\mathbb{X})$ the space of $\mathbb{X}$-valued random variables with finite $p$-th moment, i.e.,
$E\|X\|^{p}=\int_{\Omega}\|X\|^{p}\,d P<\infty.$ Then equipped it  with the norm $\|X\|_p:=(E\|X\|^{p})^{\frac{1}{p}}$, it is a Banach space. Let $\mathcal{L}^{\infty}(\mathbb{R}, \mathcal{L}^p(\Omega,\mathbb{X})$ be the space of measurable functions $g:\mathbb{R}\to \mathcal{L}^p(\Omega,\mathbb{X})$ with finite essential supremum norm
$\|g\|_{\infty}:=\operatorname*{ess\,sup}_{t\in\mathbb{R}}\|g(t)\|_{p}$,
this space is a Banach space.

\begin{definition}\label{Def2.1}
A function  $f: \mathbb{R} \times\mathcal{L}^p(\Omega, \mathbb{X}) \to \mathcal{L}^p(\Omega, \mathbb{X})$ is said to be $\mathcal{L}^p$-continuous if for any $s \in \mathbb{R}$,
$\lim_{t \to s} E\|f(t,x) - f(s,x)\|^p = 0$ uniformly in $x\in \mathcal{L}^p(\Omega, \mathbb{X})$, and it is said to be $\mathcal{L}^p$-bounded if $\|f(\cdot,x)\|_0:=\sup_{t \in \mathbb{R}} E\|f(t,x)\|_p < \infty$  uniformly in $x\in \mathcal{L}^p(\Omega, \mathbb{X})$.
\end{definition}

\begin{definition}\cite{2}
 A jointly continuous function $f:\mathbb{R}\times \mathcal{L}^p(\Omega, \mathbb{X})\rightarrow\mathcal{L}^{p}(\Omega,\mathbb{X})$ is called (Bochner) almost automorphic  if for every sequence $\{s'_{k}\}$ there exist a subsequence $\{s_{k}\}$ and another process $\tilde{f}(t,x)$ such that for each $t\in\mathbb{R}$,
\[
\lim_{k\to\infty}\|f(t+s_{k},x)-\tilde{f}(t,x)\|_{p}=0\quad\text{and}\quad \lim_{k\to\infty}\|\tilde{f}(t-s_{k},x)-f(t,x)\|_{p}=0
\]
uniformly in $x\in \mathcal{L}^p(\Omega, \mathbb{X})$. This class is $AA(\mathbb{R}\times \mathcal{L}^p(\Omega, \mathbb{X}),\mathcal{L}^{p}(\Omega,\mathbb{X}))$.
\end{definition}

\begin{remark}\label{re11}
If $f\in AA(\mathbb{R}\times \mathcal{L}^p(\Omega, \mathbb{X}),\mathcal{L}^{p}(\Omega,\mathbb{X}))$, then $\|f(\cdot,x)\|_0=\|\tilde{f}(\cdot,x)\|_0$, in addition, if there is a constant $L_f$ such that for any $x,y\in \mathcal{L}^p(\Omega, \mathbb{X})$, $\|f(\cdot,x)-f(\cdot,x)\leq L_f\|x-y\|$, then $\|\tilde{f}(\cdot,x)-\tilde{f}(\cdot,x)\leq L_f\|x-y\|$.
\end{remark}

Assume that $\mathcal{L}_{\mathrm{loc}}^{p}(\mathbb{R},\mathcal{L}^{p}(\Omega,\mathbb{X}))$ is the space of all stochastic processes $X:\mathbb{R}\to\mathcal{L}^{p}(\Omega,\mathbb{X})$ that are locally $p$-integrable.
For $X\in\mathcal{L}_{\mathrm{loc}}^{p}(\mathbb{R},\mathcal{L}^{p}(\Omega,\mathbb{X}))$, its Besicovitch seminorm is defined as
\[
\|X\|_{B^{p}}
= \left(\overline{ \lim_{l\to\infty}} \frac{1}{2l}\int_{-l}^{l} \|X(t)\|_p^{p}\,\dd t \right)^{1/p}.
\]

\begin{definition}
A process $X\in \mathcal{L}_{\mathrm{loc}}^{p}(\mathbb{R},\mathcal{L}^{p}(\Omega,\mathbb{X}))$ is called $B^{p}$-bounded if $\|X\|_{B^{p}}<\infty$. We write $B_{B}^{p}(\mathbb{R},\mathcal{L}^{p}(\Omega,\mathbb{X}))$ for the class of all such processes.
\end{definition}

\begin{definition}
A process $X\in \mathcal{L}_{\mathrm{loc}}^{p}(\mathbb{R},\mathcal{L}^{p}(\Omega,\mathbb{X}))$ is called $B^{p}$-continuous if
$\lim_{\varpi\to0}\|X(\cdot+\varpi)-X(\cdot)\|_{B^{p}}=0.$
The corresponding class is denoted by $B_{C}^{p}(\mathbb{R},\mathcal{L}^{p}(\Omega,\mathbb{X}))$.
\end{definition}

Let $B_{BC}^{p}(\mathbb{R},\mathcal{L}^{p}(\Omega,\mathbb{X}))=B_{B}^{p}(\mathbb{R},\mathcal{L}^{p}(\Omega,\mathbb{X}))\cap B_{C}^{p}(\mathbb{R},\mathcal{L}^{p}(\Omega,\mathbb{X}))$. We denote by $B_{BC}^{p}(\mathbb{R}\times \mathcal{L}^{p}(\Omega,\mathbb{X}),\mathcal{L}^{p}(\Omega,\mathbb{X}))$ the set of functions $f: \mathbb{R}\times \mathcal{L}^{p}(\Omega,\mathbb{X})\rightarrow \mathcal{L}^{p}(\Omega,\mathbb{X})$ such that $f(\cdot,x)\in B_{BC}^{p}(\mathbb{R},\mathcal{L}^{p}(\Omega,\mathbb{X}))$ for all $x\in \mathcal{L}^{p}(\Omega,\mathbb{X})$.

\begin{lemma}\label{lemm21}
Let $\mathbb{V}
:= \bigl\{ y\in B_{BC}^{p}(\mathbb{R},\mathcal{L}^{p}(\Omega,\mathbb{X}))
\cap \mathcal{L}^{\infty}(\mathbb{R},\mathcal{L}^{p}(\Omega,\mathbb{X})) \bigr\}.$
Equipped with the norm $\|\cdot\|_{\infty}$, the space $\mathbb{V}$ is a Banach space.
\end{lemma}

\begin{proof}
Since  the space $\mathcal{L}^{p}(\Omega,\mathbb{X})$ is a Banach space. Hence $\mathcal{L}^{\infty}(\mathbb{R},\mathcal{L}^{p}(\Omega,\mathbb{X}))$, equipped with the essential supremum norm $\|\cdot\|_{\infty}$, is a Banach space.

It is enough to show that $\mathbb{V}$ is closed in $\mathcal{L}^{\infty}(\mathbb{R},\mathcal{L}^{p}(\Omega,\mathbb{X}))$. Let $\{y_n\}\subset \mathbb{V}$ and suppose that $y_n\to y$ in $\mathcal{L}^{\infty}(\mathbb{R},\mathcal{L}^{p}(\Omega,\mathbb{X}))$, i.e.
$\|y_n-y\|_{\infty}\to 0$ as $n\rightarrow\infty$.
Then $y\in \mathcal{L}^{\infty}(\mathbb{R},\mathcal{L}^{p}(\Omega,\mathbb{X}))$. Also, for every compact interval $I\subset\mathbb{R}$,
$\int_I \|y(t)\|_p^p\,dt\le |I|\,\|y\|_{\infty}^p<\infty,$
so $y\in \mathcal{L}_{\mathrm{loc}}^{p}(\mathbb{R},\mathcal{L}^{p}(\Omega,\mathbb{X}))$.

Moreover,
\[
\|y-y_n\|_{B^p}\le \|y-y_n\|_{\infty}\to 0 \text{ as } n\rightarrow\infty.
\]
Thus, for some fixed $n$,
\[
\|y\|_{B^p}\le \|y_n\|_{B^p}+\|y-y_n\|_{B^p}\le \|y_n\|_{B^p}+\|y-y_n\|_{\infty}<\infty,
\]
so $y\in B_B^p(\mathbb{R},\mathcal{L}^{p}(\Omega,\mathbb{X}))$.

For $B^p$-continuity, let $\varepsilon>0$. Choose $n$ such that
\[
\|y-y_n\|_{\infty}<\frac{\varepsilon}{4}.
\]
Since $y_n\in B_C^p$, there exists $\delta>0$ such that
\[
\|y_n(\cdot+\varpi)-y_n(\cdot)\|_{B^p}<\frac{\varepsilon}{2}
\]
whenever $|\varpi|<\delta$. Then
\begin{align*}
\|y(\cdot+\varpi)-y(\cdot)\|_{B^p}
&\le \|y_n(\cdot+\varpi)-y_n(\cdot)\|_{B^p}+2\|y-y_n\|_{\infty}
<\frac{\varepsilon}{2}+2\cdot\frac{\varepsilon}{4}
=\varepsilon.
\end{align*}
Hence $y\in B_C^p(\mathbb{R},\mathcal{L}^{p}(\Omega,\mathbb{X}))$. Therefore
$y\in B_{BC}^p\cap \mathcal{L}^{\infty}(\mathbb{R},\mathcal{L}^{p}(\Omega,\mathbb{X}))=\mathbb{V}.$
So $\mathbb{V}$ is closed in the Banach space $\mathcal{L}^{\infty}(\mathbb{R},\mathcal{L}^{p}(\Omega,\mathbb{X}))$. A closed subspace of a Banach space is itself a Banach space, and therefore $(\mathbb{V},\|\cdot\|_{\infty})$ is a Banach space. This completes the proof.
\end{proof}

\begin{definition}\label{d1}\cite{a1}
A process $X\in B_{BC}^{p}(\mathbb{R},\mathcal{L}^{p}(\Omega,\mathbb{X}))$ is Besicovitch almost automorphic in $p$-th mean if for every sequence $\{s'_{k}\}$ there exist a subsequence $\{s_{k}\}$ and another process $\tilde{X}$ such that
\[
\lim_{k\to\infty}\|X(\cdot+s_{k})-\tilde{X}(\cdot)\|_{B^p}=0\quad\text{and}\quad \lim_{k\to\infty}\|\tilde{X}(\cdot-s_{k})-X(\cdot)\|_{B^p}=0.
\]
The collection is $B^pAA_m(\mathbb{R},\mathcal{L}^{p}(\Omega,\mathbb{X}))$.
\end{definition}

\begin{definition}\label{ddd1}
A function $f\in B_{BC}^{p}(\mathbb{R}\times \mathcal{L}^{p}(\Omega,\mathbb{X}),\mathcal{L}^{p}(\Omega,\mathbb{X}))$ is Besicovitch almost automorphic in $p$-th mean if for every sequence $\{s'_{k}\}$ there exist a subsequence $\{s_{k}\}$ and another process $\tilde{X}$ such that
\[
\lim_{k\to\infty}\|f(\cdot+s_{k},x)-\tilde{f}(\cdot,x)\|_{B^p}=0\quad\text{and}\quad \lim_{k\to\infty}\|\tilde{f}(\cdot-s_{k},x)-f(\cdot,x)\|_{B^p}=0
\]
uniformly in $x\in \mathbb{D}$, where $\mathbb{D}$ is any bounded subset of $\mathcal{L}^{p}(\Omega,\mathbb{X})$. The collection of such functions  is denoted by $B^pAA_m(\mathbb{R}\times \mathcal{L}^{p}(\Omega,\mathbb{X}),\mathcal{L}^{p}(\Omega,\mathbb{X}))$.
\end{definition}

\begin{remark}
Obviously, $AA(\mathbb{R}\times \mathcal{L}^{p}(\Omega,\mathbb{X}),\mathcal{L}^{p}(\Omega,\mathbb{X}))\subset B^pAA_m(\mathbb{R}\times \mathcal{L}^{p}(\Omega,\mathbb{X}),\mathcal{L}^{p}(\Omega,\mathbb{X}))$.
\end{remark}

Let $(\mathbb{E},d)$ be a separable complete metric space, $\mathcal{B}(\mathbb{E})$ its Borel $\sigma$-algebra, and $\mathcal{P}(\mathbb{E})$ the set of probability measures on $\mathcal{B}(\mathbb{E})$. Let $B_C(\mathbb{E})$ be the space of bounded continuous functions $g:\mathbb{E}\to\mathbb{R}$ with $\|g\|_{B}:=\sup_{x\in\mathbb{E}}|g(x)|$. For $g\in B_C(\mathbb{E})$ and $\zeta,\varsigma\in \mathcal{P}(\mathbb{E})$ define
\[
\|g\|_{L}= \sup_{a\neq b}\frac{|g(a)-g(b)|}{d(a,b)},\,\, \|g\|_{BL}= \max\{\|g\|_{B},\|g\|_{L}\},\,\,
d_{BL}(\zeta,\varsigma) = \sup_{\|g\|_{BL}\le 1}\Big|\int_{\mathbb{R}{E}} g\,\mathrm{d}(\zeta-\varsigma)\Big|.
\]

\begin{definition}\cite{a1}
A process $X$ is $p$-th Besicovitch almost automorphic in distribution if for every sequence $\{s'_{k}\}$ there exist a subsequence $\{s_{k}\}$ and another process $\tilde{X}$ such that
\[\lim_{k \to \infty }\bigg(\overline{\lim_{l\to \infty}} \frac{1}{2l}\int_{-l}^{l} d_{BL}^p(P \circ [X(t+s_k)]^{-1}, P \circ [\tilde X(t)]^{-1}) \dd t \bigg)^{\frac{1}{p}}=0 \]
	and
	\[\lim_{k \to \infty }\bigg(\overline{\lim_{l\to \infty}} \frac{1}{2l}\int_{-l}^{l} d_{BL}^p(P \circ [\tilde X(t-s_k)]^{-1}, P \circ [X(t)]^{-1}) \dd t\bigg)^{\frac{1}{p}} =0, \]
		where $P \circ [X(t)]^{-1}$ is the distribution of $X(t)$ under $P$.
The class is $B^pAA_d(\mathbb{R},\mathcal{L}^{p}(\Omega,\mathbb{O}^{n}))$.
\end{definition}

\begin{lemma}\cite{a3}\label{bdg} (Burkholder--Davis--Gundy inequality)
If $f\in\mathcal{L}^{2}([t_{0},T],\mathbb{R})$, $p\ge2$, and $B(t)$ is a Brownian motion, then
\[
E\bigg(\sup_{t\in\mathbb{R}}\bigg|\int_{t_{0}}^{t}f(s)\dd B(s)\bigg|\bigg)^{p}\le C_{p}E\bigg[\int_{t_{0}}^{T}|f(s)|^{2}\dd s\bigg]^{p/2},
\]
where $C_{p}=\big(\frac{p^{p+1}}{2(p-1)^{p-1}}\big)^{p/2}$.
\end{lemma}

\begin{lemma}\label{GR}\cite{27}
Let $\varphi:\mathbb{R}\to\mathbb{R}$ be continuous and satisfy
\[
0\le\varphi(t)\le a(t)+b\int_{-\infty}^{t}e^{-c(t-s)}\varphi(s)\dd s,
\]
where $a$ is locally integrable, $b,c>0$ with $c>b$, and the integral is assumed convergent. Then for any $d\in(0,c-b]$ with $\int_{-\infty}^{0}e^{ds}a(s)\dd s<\infty$,
\[
\varphi(t)\le a(t)+b\int_{-\infty}^{t}e^{-d(t-s)}a(s)\dd s.
\]
\end{lemma}

\section{Besicovitch-like  pseudo almost automorphic stochastic process}
\setcounter{equation}{0}
\indent

In this section, we introduce the concept of Besicovitch-like  pseudo almost automorphic stochastic processes.

Let $\Sigma$ be the Lebesgue $\sigma$-field on $\mathbb{R}$ and $\mathcal{M}_{p}$ the set of all positive measures $\mu$ on $\Sigma$ with $\mu(\mathbb{R})=+\infty$ and $\mu([a,b])<+\infty$ for all $a\le b$.

For $\mu,\nu\in\mathcal{M}_{p}$, the relation $\mu\sim\nu$ signifies that $\mu$ and $\nu$ are equivalent in the following sense: there exist constants $\alpha,\beta>0$ and a bounded interval $A$ (possibly empty) such that, for every $C\in\Sigma$ with $C\cap A=\emptyset$, the value $\mu(C)$ is sandwiched between $\alpha\nu(C)$ and $\beta\nu(C)$, namely,
$\alpha\,\nu(C)\le \mu(C)\le \beta\,\nu(C).$

For $\mu\in\mathcal{M}_{p}$ and $\sigma\in\mathbb{R}$, the shifted measure $\mu_{\sigma}$ is defined by $\mu_{\sigma}(A)=\mu(\{a+\sigma:a\in A\})$.

We impose the following two basic requirements:
\begin{itemize}
\item[$(A_1)$] $\mu,\nu\in\mathcal{M}_{p}$ satisfy $\displaystyle\overline{\lim_{h\to \infty}} \frac{\mu([-h,h])}{\nu([-h,h])}<\infty$.
\item[$(A_2)$] For every $\mu\in\mathcal{M}_{p}$ and $\eta\in\mathbb{R}$ there exist $\rho>0$ and a bounded interval $\mathcal{J}\subset\mathbb{R}$ such that
\[
\mu(\{a+\eta:a\in A\})\le\rho\,\mu(A),
\]
whenever $A\in\Sigma$ satisfies $A\cap\mathcal{J}=\emptyset$.
\end{itemize}
\begin{lemma}\label{a1}\cite{b6}
$\mu\in\mathcal{M}_{p}$ satisfies $(A_2)$ if and only if $\mu\sim\mu_{\sigma}$ for all $\sigma\in\mathbb{R}$.
\end{lemma}
\begin{lemma}\cite{b6}\label{2.5}
Condition $(A_2)$ implies that, for all $\eta > 0$,
\[
\overline{\lim_{h\to \infty}} \frac{\mu([-h - \eta, h + \eta])}{\mu([-h, h])} < \infty.
\]
\end{lemma}
\begin{lemma}\label{lem27}
Let $\mu \in \mathcal{M}_{p}$ satisfy condition $(A_2)$. Then, for every $\sigma \in \mathbb{R}$ and $\eta > 0$,
\[
\overline{\lim_{h\to \infty}} \frac{\mu_\sigma([-h-\eta,\, h+\eta])}{\mu([-h,\, h])} < \infty .
\]
\end{lemma}

\begin{proof}
By $(A_2)$ and Lemma \ref{a1}, $\mu$ is equivalent to its translate $\mu_\sigma$. Hence, there exist a constant $\rho>0$ and a bounded interval $J$ such that, for every $A \in \Sigma$ with $A \cap J = \emptyset$,
$\mu_\sigma(A) \le \rho \mu(A).$

Fix $\eta > 0$. For $h$ large enough so that $J \subset [-h-\eta, h+\eta]$, we have
\begin{align*}
 \mu_\sigma([-h-\eta,\, h+\eta])
&\le \mu_\sigma(J) + \rho\mu([-h-\eta, h+\eta]\setminus J)\\
&\le \mu_\sigma(J) + \rho\mu([-h-\eta, h+\eta]).
\end{align*}
Since $\mu_\sigma(J) = \mu(J-\sigma)$ is finite and $\mu(\mathbb{R})=+\infty$, it follows that
\[
\overline{\lim_{h\to \infty}}\frac{\mu_\sigma([-h-\eta,\, h+\eta])}{\mu([-h, h])}
\le \rho\overline{\lim_{h\to \infty}} \frac{\mu([-h-\eta,\, h+\eta])}{\mu([-h, h])},
\]
which, combined with Lemma \ref{2.5}, implies that the desired upper limit is finite. This completes the proof.
\end{proof}

\begin{definition}\label{ddf1}
Let $\mu,\nu\in\mathcal{M}_{p}$. A process  $X: \mathbb{R} \to \mathcal{L}^{p}(\Omega,\mathbb{X})$ belongs to $\varepsilon_M(\mathbb{R},\mathcal{L}^{p}(\Omega,\mathbb{X}),
\mu,\nu)$ if and only if
\[
\lim_{h\to\infty}\frac{1}{\nu([-h,h])}\int_{-h}^{h}E\|X(t)\|^{p}\dd\mu(t)=0.
\]
\end{definition}
\begin{definition}\label{ddf2}
Let $\mu,\nu\in\mathcal{M}_{p}$. A process $X:\mathbb{R}\to\mathcal{L}^{p}(\Omega,\mathbb{X})$ belongs to $\varepsilon_D(\mathbb{R},\mathcal{L}^{p}(\Omega,\mathbb{X}),\mu,\nu)$ if and only if
\[
\lim_{h\to\infty}\frac{1}{\nu([-h,h])}\int_{-h}^{h} d_{BL}^{p}\big(P\circ[X(t)]^{-1},\delta_{0}\big)\dd\mu(t)=0,
\]
 where $\delta _{0}$ denotes the distribution of random variable $\mathbf{0}$.
\end{definition}

Since
\begin{align*}
	\frac{1}{\nu ([-h,h])}\int_{-h}^{h}d_{BL}^p(P \circ [X(t)]^{-1},\delta_0)\dd \mu (t)
	\le &\frac{1}{\nu ([-h,h])}\int_{-h}^{h} \bigg(\int_\Omega  \|X(t)\|\dd P \bigg)^p\dd \mu (t)\nonumber\\
	\leq&\frac{1}{\nu ([-h,h])}\int_{-h}^{h}  E\|X(t)\|^p \dd \mu (t),
\end{align*}
from Definitions \ref{ddf1} and \ref{ddf2}, we have
\begin{theorem}\label{key}
  $\varepsilon_M(\mathbb{R},\mathcal{L}^{p}(\Omega,\mathbb{X}),\mu,\nu)\subset \varepsilon_D(\mathbb{R},\mathcal{L}^{p}(\Omega,\mathbb{X}),
\mu,\nu)$.
\end{theorem}

\begin{definition}
Let $\mu,\nu\in\mathcal{M}_{p}$. A function  $f: \mathbb{R}\times \mathcal{L}^{p}(\Omega,\mathbb{X}) \to \mathcal{L}^{p}(\Omega,\mathbb{X})$ belongs to $\varepsilon_M(\mathbb{R}\times \mathcal{L}^{p}(\Omega,\mathbb{X}),\mathcal{L}^{p}(\Omega,\mathbb{X}),
\mu,\nu)$ if and only if
\[
\lim_{h\to\infty}\frac{1}{\nu([-h,h])}\int_{-h}^{h}E\|f(t,x)\|^{p}\dd\mu(t)=0
\]
uniformly in $x\in  \mathcal{L}^{p}(\Omega,\mathbb{X})$.
\end{definition}

\begin{definition}
Let $\mu,\nu\in\mathcal{M}_{p}$. A function $f:\mathbb{R}\times \mathcal{L}^{p}(\Omega,\mathbb{X}) \to\mathcal{L}^{p}(\Omega,\mathbb{X})$ belongs to $\varepsilon_D(\mathbb{R}\times\mathcal{L}^{p}(\Omega,\mathbb{X}),\mathcal{L}^{p}(\Omega,\mathbb{X}),\mu,\nu)$ if and only if
\[
\lim_{h\to\infty}\frac{1}{\nu([-h,h])}\int_{-h}^{h} d_{BL}^{p}\big(P\circ[f(t,x)]^{-1},\delta_{0}\big)\dd\mu(t)=0
\]
uniformly in $x\in  \mathcal{L}^{p}(\Omega,\mathbb{X})$,  and $\delta _{0}$ denotes the distribution of random variable $\mathbf{0}$.
\end{definition}

Now we define the new classes of processes and functions.
\begin{definition}\label{d2}
Let $\mu,\nu\in\mathcal{M}_{p}$. A process $X:\mathbb{R}\to\mathcal{L}^{p}(\Omega,\mathbb{X})$ is called  Besicovitch-like $(\mu,\nu)$-pseudo almost automorphic in $p$-th mean if it can be decomposed as
$X = X_{1}+X_{2},$
where $X_{1}\in B^pAA_m(\mathbb{R},\mathcal{L}^{p}(\Omega,\mathbb{X}))$ and $X_{2}\in\varepsilon_M(\mathbb{R},\mathcal{L}^{p}(\Omega,\mathbb{X}),\mu,\nu)$. This class is denoted by $PB^{p}AA_m(\mathbb{R},\mathcal{L}^{p}(\Omega,\mathbb{X}),\mu,\nu)$.
\end{definition}
\begin{definition}
Let $\mu,\nu\in\mathcal{M}_{p}$. A process $X:\mathbb{R}\to\mathcal{L}^{p}(\Omega,\mathbb{X})$ is said to be Besicovitch-like $(\mu,\nu)$-pseudo almost automorphic in distribution if it can be written as
$X = X_{1}+X_{2}$
with $X_{1}\in B^pAA_d(\mathbb{R},\mathcal{L}^{p}(\Omega,\mathbb{X}))$ and $X_{2}\in\varepsilon_D(\mathbb{R},\mathcal{L}^{p}(\Omega,\mathbb{X}),\mu,\nu)$. The family of such processes is indicated by $PB^pAA_d(\mathbb{R},
\mathcal{L}^{p}(\Omega,\mathbb{X}),\mu,\nu)$.
\end{definition}

\begin{definition}\label{d2}
Let $\mu,\nu\in\mathcal{M}_{p}$. A function $f:\mathbb{R}\times\mathcal{L}^{p}(\Omega,\mathbb{X}) \to\mathcal{L}^{p}(\Omega,\mathbb{X})$ is called  Besicovitch-like $(\mu,\nu)$-pseudo almost automorphic in $p$-th mean if it can be decomposed as
$f = f_{1}+f_{2},$
where $f_{1}\in B^pAA_m(\mathbb{R}\times\mathcal{L}^{p}(\Omega,\mathbb{X}),\mathcal{L}^{p}(\Omega,\mathbb{X}))$ and $f_{2}\in\varepsilon_M(\mathbb{R}\times\mathcal{L}^{p}(\Omega,\mathbb{X}),\mathcal{L}^{p}(\Omega,\mathbb{X}),\mu,\nu)$. This class is denoted by $B^pAA_m(\mathbb{R}\times\mathcal{L}^{p}(\Omega,\mathbb{X}),\mathcal{L}^{p}(\Omega,\mathbb{X}),\mu,\nu)$.
\end{definition}
\begin{definition}
Let $\mu,\nu\in\mathcal{M}_{p}$. A function $f:\mathbb{R}\times\mathcal{L}^{p}(\Omega,\mathbb{X})\to\mathcal{L}^{p}(\Omega,\mathbb{X})$ is called Besicovitch-like $(\mu,\nu)$-pseudo almost automorphic in distribution if it can be written as
$f = f_{1}+f_{2}$
with $f_{1}\in B^pAA_d(\mathbb{R}\times\mathcal{L}^{p}(\Omega,\mathbb{X}),\mathcal{L}^{p}(\Omega,\mathbb{X}))$ and $f_{2}\in\varepsilon_D(\mathbb{R}\times\mathcal{L}^{p}(\Omega,\mathbb{X}),\mathcal{L}^{p}(\Omega,\mathbb{X}),\mu,\nu)$. The collection is $DB_{PAA}^{p}(\mathbb{R}\times\mathcal{L}^{p}(\Omega,\mathbb{X}),
\mathcal{L}^{p}(\Omega,\mathbb{X}),\mu,\nu)$.
\end{definition}

\section{Besicovitch-like  pseudo almost automorphic solutions in distribution}
\setcounter{equation}{0}
\indent	

In this section, we will establish the existence, uniqueness and global exponential stability of Besicovitch-like $(\mu,\nu)$-pseudo almost automorphic solutions in distribution to \eqref{1}.

\begin{definition}\label{d1}
An $\mathcal F_t$-progressively measurable stochastic process $X(t)$ is called a solution of system \eqref{1},  if $X(t)$ satisfies the following stochastic integral equation:
\begin{align}\label{jf1}
 X(t)=T(t-t_0)X(t_0)+\int_{t_0}^tT(t-s)f(s,X(s))\dd s+\int_{t_0}^t T(t-s)g(s,X(s))\dd W(s),
\end{align}
for all $t\ge t_0$ and each $t_0\in\mathbb R.$
\end{definition}

Consider the following space:
\[
\mathbb{U}:= \big\{ y\in B_{BC}^{p}(\mathbb{R},\mathcal{L}^{p}(\Omega,\mathbb{X})) \cap \mathcal{L}^{\infty}(\mathbb{R},\mathcal{L}^{p}(\Omega,\mathbb{X})) \text{ is } \mathcal{F}_t-\text{progressively measurable}  \big\}
\]
equipped with the norm $\|y\|_{\infty}=\operatorname*{ess\,sup}_{t\in\mathbb{R}}\|y(t)\|_{p}$. Then, by Lemma \ref{lemm21}, it is easy to see that $\mathbb{U}$ is a Banach space.

We  make the following assumptions:
\begin{itemize}
  \item [$(H_{1})$]
 The operator  $A :Dom(A)\subset \mathbb X\rightarrow \mathbb X$  generates an exponentially stable semigroup $\{T(t); t\in \mathbb{R}^+\},$ i.e., there are constants $M\geq 1$ and $\delta>0$ such that
   $\|T(t) \|\le Me^{-\delta t}$ for all $t\in \mathbb{R}^+$.

   \item [$(H_{2})$]
Functions $f\in B^pAA_m(\mathbb{R}\times \mathcal{L}^{p}(\Omega,\mathbb{X}),\mathcal{L}^{p}(\Omega,\mathbb{X}),\mu,\nu)$ and $g\in B^pAA_m(\mathbb{R}\times \mathcal{L}^{p}(\Omega,\mathbb{X}),\mathbb{L}_2^0,\mu,\nu)$ satisfy that for all  $x,y\in  \mathcal{L}^{p}(\Omega,\mathbb{X}))$ and $t\in\mathbb R$, there exist constants $K_f,K_g,L_f,L_g$ such that
\begin{gather*}
  E\| f(t,x)\|^p\leq K_f(1+E\| x\|^p), \,\,
  E \| g(t,x)\|^p_{\mathbb{L}_2^0}\leq K_g (1+E\| x\|^p),\\
  E\| f(t,x)-f(t,y)\|\leq L_fE\| x-y\|, \,\,
  E \| g(t,x)-g(t,y)\|_{\mathbb{L}_2^0}\leq L_gE\| x-y\|,\\
   E\| \tilde{f}(t,x)\|^p\leq K_f(1+E\| x\|^p), \,\,
  E \| \tilde{g}(t,x)\|^p_{\mathbb{L}_2^0}\leq K_g (1+E\| x\|^p),\\
  E\| \tilde{f}(t,x)-\tilde{f}(t,y)\|\leq L_fE\| x-y\|, \,\,
  E \| \tilde{g}(t,x)-\tilde{g}(t,y)\|_{\mathbb{L}_2^0}\leq L_gE\| x-y\|,
\end{gather*}
where $\tilde{f},\tilde{g}$ are defined similarly as in Definition \ref{ddd1}.

 \item [$(H_{3})$]
Functions $f\in PB^pAA_m(\mathbb{R}\times \mathcal{L}^{p}(\Omega,\mathbb{X}),\mathcal{L}^{p}(\Omega,\mathbb{X}),\mu,\nu)$ and $g\in PB^pAA_m(\mathbb{R}\times \mathcal{L}^{p}(\Omega,\mathbb{X}),\mathbb{L}_2^0,\\
\mu,\nu)$
can be decomposed as
\[
f=f_1+f_2 \text{ and } g=g_1+g_2,
\]
where $f_1\in B^pAA_m(\mathbb{R}\times \mathcal{L}^{p}(\Omega,\mathbb{X}),\mathcal{L}^{p}(\Omega,\mathbb{X}))$, $g_1\in B^pAA_m(\mathbb{R}\times \mathcal{L}^{p}(\Omega,\mathbb{X}),\mathbb{L}_2^0)$, $f_2\in\varepsilon^{p}_M(\mathbb{R}\times \mathcal{L}^{p}(\Omega,\mathbb{X}),\mathcal{L}^{p}(\Omega,\mathbb{X}),\mu,\nu) $, and $g_2\in \varepsilon_M(\mathbb{R}\times \mathcal{L}^{p}(\Omega,\mathbb{X}),\mathbb{L}_2^0,\mu,\nu).$

\item [$(H_{4})$] Functions $f$ and $g$ satisfy $(H_3)$, and
   for all  $x,y\in  \mathcal{L}^{p}(\Omega,\mathbb{X}))$ and $t\in\mathbb R$, there exist constants $K_f,K_g,L_f,L_g$ such that
\begin{gather*}
E\| f_1(t,x)\|^p\leq K_f(1+E\| x\|^p), \,\,
E \| g_1(t,x)\|^p_{\mathbb{L}_2^0}\leq K_g (1+E\| x\|^p),\\
E\| f_1(t,x)-f_1(t,y)\|\leq L_fE\| x-y\|, \,\, E \| g_1(t,x)-g_1(t,y)\|_{\mathbb{L}_2^0}\leq L_gE\| x-y\|,\\
E\| \tilde{f}_1(t,x)\|^p\leq K_f(1+E\| x\|^p), \,\,
  E \| \tilde{g}_1(t,x)\|^p_{\mathbb{L}_2^0}\leq K_g (1+E\| x\|^p),\\
  E\| \tilde{f}_1(t,x)-\tilde{f}_1(t,y)\|\leq L_fE\| x-y\|, \,\,
  E \| \tilde{g}_1(t,x)-\tilde{g}_1(t,y)\|_{\mathbb{L}_2^0}\leq L_gE\| x-y\|,
\end{gather*}
where $\tilde{f}_1,\tilde{g}_1$ are defined similarly as in Definition \ref{ddd1}.
\item [$(H_{5})$]Let
 \begin{align*}
  R_1 :=&4^{p-1}M^p\bigg[ L_f^p\bigg(\frac{1}{\delta}\bigg)^{p}+  C_p  L_g^p\bigg(\frac{p-2}{p\delta}\bigg)^{\frac{p-2}{2}}\bigg(\frac{2}{p\delta}\bigg)\bigg]<1,\,\, (p>2),\\
 R_2:=&4 M^2 \bigg[L_f^2\bigg(\frac{1}{\delta}\bigg)^{2}+ L_g^2 \bigg(\frac{1}{2\delta}\bigg)\bigg]<1.
\end{align*}

\item [$(H_{6})$]Let
  \begin{align*}
  \beta_1:=&4^{p-1}M^p\bigg[L_f^p\bigg(\frac{1}{\delta}\bigg)^{p-1}+ C_p   L_g^p\bigg(\frac{p-2}{p\delta}\bigg)^{\frac{p-2}{2}}\bigg]<\delta, \,\, (p>2),\\
  \beta_2:=&4 M^2\bigg[L_f^2\bigg(\frac{1}{\delta}\bigg) + L_g^2\bigg]<\delta.
\end{align*}
\end{itemize}
 \begin{remark}
 Clearly, condition $(H_6)$ implies condition $(H_5)$.
 \end{remark}

\begin{lemma}\label{th3.1}
Under the conditions $(H_1)$, $(H_2)$ and $(H_5)$,  system \eqref{1} admits a unique solution in $\mathbb{U}$.
\end{lemma}
\begin{proof}
 By  \eqref{jf1} and condition $(H_1)$, letting $t_0\rightarrow -\infty$ yields the following stochastic integral equation
\begin{align}\label{eq1}
X(t)=\int_{-\infty}^t T(t-s)f(s,X(s))ds+\int_{-\infty}^t T(t-s)g(s,X(s))\mathrm{d} W(s),
\end{align}
which is also a solution of \eqref{1}.

We define an operator $\Psi:\mathbb X\rightarrow\mathbb X$ by setting
\begin{align*}
(\Psi \phi)(t)=&\int_{-\infty}^t T(t-s)f(s,\phi(s))\dd s+\int_{-\infty}^t T(t-s)g(s,\phi(s))\dd W(s),
\end{align*}
where $\phi\in \mathbb{X}, t\in \mathbb{R}.$

Firstly, we will show that $\Psi$ maps $\mathbb{U}$ into $\mathbb{U}$. For any $\varphi \in \mathbb{U}$, we have
\begin{align}\label{6}
 E\| (\Psi\phi)(t)\|^p
\leq& 2^{p-1}E\bigg\|\int_{-\infty}^t T(t-s)f(s,\phi(s))\mathrm{d}s \bigg\|^p
+ 2^{p-1}E\bigg\|\int_{-\infty}^t T(t-s)g(s,\phi(s))\mathrm{d}W(s) \bigg\|^p\nonumber\\
:=&G_{1}+G_{2} .
\end{align}
By the H\"{o}lder inequality,  $(H_1)$ and $(H_2)$, one derives that
\begin{align}\label{8}
G_{1}
\le& 2^{p-1}\bigg(\int_{-\infty }^{t}\|T(t-s)\| \mathrm{d}s\bigg)^{p-1}E\int_{-\infty }^{t}\|T(t-s)\| \| f(s,\phi(s))\|^{p}\mathrm{d}s \nonumber\\
\le&2^{p-1}M^{\frac{p +q}{q}} \bigg(\int_{-\infty }^{t} e^{-\delta(t-s)}\mathrm{d}s\bigg)^{p-1}\int_{-\infty }^{t}e^{-\delta(t-s)} K_f^p E(1+\|\phi(s)\|^{p})\mathrm{d}s \nonumber\\
\le&2^{p-1}M^{p} K_f^p \bigg(\frac{1}{\delta}\bigg)^{p} (1+\|\phi\|_{\infty}^p) <\infty.
\end{align}

 When $p>2$, by Lemma \ref{bdg}, conditions $(H_1)$ and $(H_2)$, we can deduce that
\begin{align}\label{5}
G_{2}
\le&2^{p-1}C_pE\bigg(\int_{-\infty}^t \|T(t-s)g(s,\phi(s))\|^2\mathrm{d}s \bigg)^{\frac{p}{2}}\nonumber\\
\le&2^{p-1}C_p \bigg(\int_{-\infty}^{t} \|T(t-s) \|^{\frac{p}{p-2}} \mathrm{d}s\bigg)^{\frac{p-2}{p}\times\frac{p}{2}}
E \int_{-\infty}^{t} \|T(t-s) \|^{\frac{p}{2}} \|g(s,\phi(s))\|_{\mathbb{L}_2^0}^p \mathrm{d}s \nonumber\\
\le&2^{p-1}C_pM^{p} K_g^p \bigg(\frac{p-2}{p\delta}\bigg)^{\frac{p-2}{2}}
\frac{2}{p\delta}  (1+\|\phi\|_{\infty}^p) <\infty,
\end{align}
and when $p=2$,  according to the It\^{o} isometry, we infer that
\begin{align}\label{4}
G_{2}
=&2 \int_{-\infty}^t \|T(t-s)\|^2 E  \|g(s,\phi(s)) \|_{\mathbb{L}_2^0}^2\mathrm{d}s   \nonumber\\
 \le&2M^2\int_{-\infty}^t  e^{-2\delta(t-s)} K_g^2E(1+\|\phi(s) \|^2)\mathrm{d}s\nonumber\\
 \le&2M^2K_g^2\frac{1}{2\delta}(1+\|\phi\|_{\infty}^2)<\infty.
\end{align}
Substituting \eqref{8}--\eqref{4} into \eqref{6}, we obtain
$\| \Psi\phi \|_{\infty}=\operatorname*{ess\,sup}_{t\in\mathbb{R}}\|(\Psi\phi)(t)\|_{p}<\infty,$
which implies that $\Psi$ is bounded.  Since $\|\Psi\phi\|_{B^p} \leq \|\Psi\phi\|_\infty$, $\Psi\phi$ is $B^p$- bounded.

Next, we will show that $\Psi \varphi$ is $B^p$-continuous. For any $\phi \in \mathbb{U}^*$ and $\varpi>0$, we have
\begin{align}\label{c0}
&\|(\Psi\phi)(\cdot+\varpi)-(\Psi\phi)(\cdot)\|_{B^p}^p \nonumber\\
\leq&4^{p-1}\overline{\lim_{l\to \infty}}\frac{1}{2l}\int_{-l}^{l}E\bigg\|\int_{-\infty}^{t} [T(t+\varpi-s)-T(t-s)]f(s,\phi(s))\dd s\bigg\|^p\dd t\nonumber\\
&+4^{p-1}\overline{\lim_{l\to \infty}}\frac{1}{2l}\int_{-l}^lE\bigg\|\int_{t}^{t+\varpi} T(t+\varpi-s)f(s,\phi(s))\dd s\bigg\|^p\dd t\nonumber\\
&+4^{p-1}\overline{\lim_{l\to \infty}}\frac{1}{2l}\int_{-l}^{l}E\bigg\|\int_{-\infty}^{t}[ T(t+\varpi-s)-T(t-s)]g(s,\phi(s))\dd W(s)\bigg\|^p\dd t\nonumber\\
&+4^{p-1}\overline{\lim_{l\to \infty}}\frac{1}{2l}\int_{-l}^{l}E\bigg\|\int_{t}^{t+\varpi} T(t+\varpi-s)g(s,\phi(s))\dd W(s) \bigg\|^p\dd t\nonumber\\
=:&A_1+A_2+A_3+A_4.
\end{align}

According to the H\"older inequality and the Fubini theorem, we can infer that
\begin{align}
  A_1  \leq&4^{p-1}M^p\overline{\lim_{l\to \infty}}\frac{1}{2l}\int_{-l}^{l}E\bigg(\int_{-\infty}^{t}e^{-\delta(t-s)} \|T(\varpi)-T(0)\|\|f(s,\phi(s))\|\dd s\bigg)^p\dd t\nonumber\\
  \leq&4^{p-1}M^p\overline{\lim_{l\to \infty}}\frac{1}{2l}\int_{-l}^{l}E\bigg(\int_{-\infty}^{t}e^{-\delta(t-s)}\dd s\bigg)^{p-1} \int_{-\infty}^{t}e^{-\delta(t-s)}\|T(\varpi)-T(0)\|^p\|f(s,\phi(s))\|^p\dd s \dd t\nonumber\\
  \leq&4^{p-1}M^pK_f^p\bigg(\frac{1}{\delta}\bigg)^{p-1}\|T(\varpi)-T(0)\|^p\overline{\lim_{l\to \infty}}\frac{1}{2l}\int_{-l}^{l} \int_{-\infty}^{t}e^{-\delta(t-s)}(1+E\|\phi(s)\|^p)\dd s \dd t\nonumber\\
  \leq&4^{p-1}M^pK_f^p\bigg(\frac{1}{\delta}\bigg)^{p-1}\|T(\varpi)-T(0)\|^p \int_{0}^{\infty}e^{-\delta s}\dd s\overline{\lim_{l\to \infty}}\frac{2(l+|s|)}{2l}\nonumber\\
  &\times\frac{1}{2(l+|s|)}\int_{-l-|s|}^{l+|s|}(1+E\|\phi(t)\|^p) \dd t\nonumber\\
  =&4^{p-1}M^pK_f^p\bigg(\frac{1}{\delta}\bigg)^{p}\|T(\varpi)-T(0)\|^p (1+\|\phi\|_{B^p}^p),\label{c1}\\
A_2\leq&4^{p-1}M^p\overline{\lim_{l\to \infty}}\frac{1}{2l}\int_{-l}^lE\bigg(\int_{t}^{t+\varpi}e^{-\delta(t+\varpi-s)} \|f(s,\phi(s))\|\dd s\bigg)^p\dd t\nonumber\\
\leq&4^{p-1}M^pK_f^p\overline{\lim_{l\to \infty}}\frac{1}{2l}\int_{-l}^l\bigg(\int_{t}^{t+\varpi}e^{-\delta(t+\varpi-s)}\dd s\bigg)^{p-1}\int_{t}^{t+\varpi}e^{-\delta(t+\varpi-s)}(1+E\|\phi(s)\|^p)\dd s \dd t\nonumber\\
\leq&4^{p-1}M^pK_f^p\bigg(\frac{1}{\delta}(1-e^{-\delta\varpi})\bigg)^{p}(1+E\|\phi\|_\infty^p).\label{c2}
\end{align}

When $p>2$, according to Lemma \ref{bdg} and the Fubini theorem, we have that
\begin{align}
 A_3 \leq &4^{p-1}C_p\overline{\lim_{l\to \infty}}\frac{1}{2l}\int_{-l}^{l}E\bigg(\int_{-\infty}^{t}\|T(t+\varpi-s)-T(t-s)\|^2\|g(s,\phi(s))\|^2_{\mathbb{L}_2^0}\dd s\bigg)^{\frac{p}{2}}\dd t\nonumber\\
 \leq &4^{p-1}C_p M^p\|T(\varpi)-T(0)\|^p\overline{\lim_{l\to \infty}}\frac{1}{2l}\int_{-l}^{l}E\bigg(\int_{-\infty}^{t}e^{-2\delta(t-s)}\|g(s,\phi(s))\|^2_{\mathbb{L}_2^0}\dd s\bigg)^{\frac{p}{2}}\dd t\nonumber\\
 \leq &4^{p-1}C_p M^p K^p_g\|T(\varpi)-T(0)\|^p\overline{\lim_{l\to \infty}}\frac{1}{2l}\int_{-l}^{l}\bigg(\int_{-\infty}^{t}e^{-\frac{p}{p-2}\delta(t-s)}\dd s\bigg)^{\frac{p-2}{2}}  \nonumber\\
 &\times\int_{-\infty}^{t}e^{-\frac{p}{2}\delta(t-s)}(1+ E\|\phi(s)\|^p)\dd s \dd t\nonumber\\
 \leq &4^{p-1}C_p M^p K^p_g\|T(\varpi)-T(0)\|^p\bigg(\frac{p-2}{p\delta}\bigg)^{\frac{p-2}{2}}\overline{\lim_{l\to \infty}}\frac{1}{2l}\int_{-l}^{l}  \int_{-\infty}^{t}e^{-\frac{p}{2}\delta(t-s)}(1+ E\|\phi(s)\|^p)\dd s \dd t\nonumber\\
 \leq&4^{p-1}C_p M^p K^p_g\|T(\varpi)-T(0)\|^p\bigg(\frac{p-2}{p\delta}\bigg)^{\frac{p-2}{2}} \frac{2}{p\delta}  (1+ E\|\phi\|_{B^p}^p),\label{c3}\\
A_4 \leq& 4^{p-1}M^pK_g^p\bigg(\frac{p-2}{p\delta}(1-e^{-\delta\varpi})\bigg)^{\frac{p-2}{2}}\overline{\lim_{l\to \infty}}\frac{1}{2l}\int_{-l}^{l} \int_{t}^{t+\varpi} e^{-\frac{p}{2}\delta(t+\varpi-s)}(1+E\|\phi(s)\|^p)\dd s \dd t\nonumber\\
\leq& 4^{p-1}M^pK_g^p\bigg(\frac{p-2}{p\delta}(1-e^{-\delta\varpi})\bigg)^{\frac{p-2}{2}}\bigg(\frac{2}{p\delta}\bigg)(1-e^{-\delta\varpi})  (1+\|\phi\|_{B^p}^p).\label{c4}
\end{align}

When $p=2$, by the It\^{o} isometry and the Fubini theorem, we can deduce that
\begin{align}
 A_3 \leq& 4 M^2\|T(\varpi)-T(0)\|^2\overline{\lim_{l\to \infty}}\frac{1}{2l}\int_{-l}^{l} \int_{-\infty}^{t} e^{-2\delta(t-s)}E\|g(s,\phi(s))\|^2_{\mathbb{L}_2^0}\dd s \dd t\nonumber\\
 \leq& 4 M^2K_g^2\|T(\varpi)-T(0)\|^2\bigg(\frac{1}{2\delta}\bigg) (1+E\|\phi\|^2_{B^p}),\label{c5}\\
 A_4 \leq &4 M^2K_g^2\overline{\lim_{l\to \infty}}\frac{1}{2l}\int_{-l}^{l}  \int_{t}^{t+\varpi} e^{-2\delta(t+\varpi-s)}(1+E\|\phi(s)\|^2)\dd s\dd t\nonumber\\
  \leq &4 M^2K_g^2\bigg(\frac{1}{2\delta}\bigg)(1- e^{-2\delta\varpi})(1+E\|\phi\|_\infty^2).\label{c6}
\end{align}

 For $p>2$, from \eqref{c0}-\eqref{c4}, and the continuity of $T(t)$ at $t=0$,   we have
\begin{align}
\lim_{\varpi\to0}\|(\Psi\phi)(\cdot+\varpi)-(\Psi\phi)(\cdot)\|_{B^p}^p=0.
\end{align}
For $p=q=2$, from \eqref{c0}-\eqref{c2}, \eqref{c5}, \eqref{c6},  and the continuity of $T(t)$ at $t=0$, we have
\begin{align}
\lim_{\varpi\to0}\|(\Psi\phi)(\cdot+\varpi)-(\Psi\phi)(\cdot)\|_{B^p}^2	=0.
\end{align}

Hence,   $\Psi\phi$ is $B^p$-continuous. Therefore, we conclude that $\Psi(\mathbb{U}) \subset \mathbb{U}$.

Finally, we will prove that $\Psi$ is a contraction mapping. For any $\varphi$, $\psi \in \mathbb{U}$, similar to the estimation processes above, for $p>2$,  we deduce that
\begin{align*}
  \|(\Psi\varphi)(t)-(\Psi\psi)(t)\|_p^p
  \leq&2^{p-1}M^p L^p_f  \bigg(\frac{1}{\delta}\bigg)^{p-1}\int_{-\infty}^t e^{-\delta(t-s)}E\|\phi(s)-\psi(s)\|^p\dd s\nonumber\\
  &+2^{p-1}C_pM^p L_g^p \bigg(\frac{p-2}{p\delta}\bigg)^{\frac{p-2}{2}}\int_{-\infty}^t e^{-\frac{p}{2}\delta(t-s)}E\|\phi(s)-\psi(s)\|^p\dd s \nonumber\\
  \leq&\frac{1}{2^{p-1}}R_1 \|\phi-\psi\|_\infty^p,
\end{align*}
where $R_1$ is defined in $(H_5)$, and for $p=2$, we can get
\begin{align*}
  &\|(\Psi\varphi)(t)-(\Psi\psi)(t)\|_p^2  \nonumber\\
  \leq&2M^2 L^2_f  \bigg(\frac{1}{\delta}\bigg) \int_{-\infty}^t e^{-\delta(t-s)}E\|\phi(s)-\psi(s)\|^2\dd s +2 M^2L_g^2  \int_{-\infty}^t e^{-2\delta(t-s)}E\|\phi(s)-\psi(s)\|^2\dd s \nonumber\\
  \leq&\frac{1}{2}R_2  \|\phi-\psi\|_\infty^2.
\end{align*}
in which $R_2$ is defined in $(H_5)$

 Consequently, for $p \ge 2$,  by $(H_5)$, we see  that $\Psi$ is a contraction mapping. Thus, $\Psi$ has a unique fixed point $x$ in $\mathbb{U}$. Subsequently, \eqref{1} admits a unique   solution $x \in \mathbb{U}$. The proof is completed.
\end{proof}

\begin{theorem}\label{th3.2}
	If conditions $(H_1)$, $(H_2)$ and $(H_5)$  are satisfied, then system \eqref{1} has a unique $p$-th Besicovitch almost automorphic solution in distribution.
\end{theorem}
\begin{proof}
According to Lemma \ref{th3.1}, system $\eqref{1}$
has a unique solution $Y(t)$ in $\mathbb{U}$ satisfying \eqref{eq1}. To finish the proof, it suffices to show that
 $Y(t)$ is $p$-th Besicovitch almost automorphic in distributions.
	
 According to $(H_2)$, let $(\alpha''_k)$ be an arbitrary sequence of real numbers, then there exists a subsequence $(\alpha'_k)$ of $(\alpha''_k)$ such that
\begin{align}
 & \lim_{k \to \infty } \| f(\cdot+\alpha'_k,x)-\tilde f(\cdot,x) \|_{B^p}^p = 0,\quad \lim_{k \to \infty } \| g(\cdot+\alpha'_k,x)-\tilde g(\cdot,x) \|_{B^p}^p = 0,\label{jh1}\\
 & \lim_{k \to \infty } \| \tilde f(\cdot-\alpha'_k,x)-f(\cdot,x) \|_{B^p}^p = 0, \quad \lim_{k \to \infty } \| \tilde g(\cdot-\alpha'_k,x) - g(\cdot,x) \|_{B^p}^p = 0,\label{jh2}
\end{align}
uniformly in $x\in \mathbb{D}$, where $\mathbb{D}$ is any bounded subset of $\mathcal{L}^p(\Omega,\mathbb{X})$.

Consider the following integral system:
\begin{align}\label{ef2}
Y(t)=\int_{-\infty}^t T(t-s)\tilde{f}(s,Y(s))ds+\int_{-\infty}^t T(t-s)\tilde{g}(s,Y(s))\mathrm{d} W(s).
\end{align}

In view of Lemma \ref{th3.1}, one can see that \eqref{ef2} has a unique solution $\tilde Y\in \mathbb{U}^*(\mathbb{R},\mathcal{L}^p(\Omega,\mathbb{O}^n))$.

From \eqref{eq1}, we can get
\begin{align*}
Y(t+\alpha_k)=&\int_{-\infty}^{t} T(t-s)f(s+\alpha_k,Y(s+\alpha_k))ds\\
&+\int_{-\infty}^{t} T(t-s)g(s+\alpha_k,Y(s+\alpha_k))\mathrm{d} [W(s+\alpha_k)-W(\alpha_k)],
\end{align*}
  where $W(s+\alpha_k)-W(\alpha_k)$ is a Brownian motion with the same distribution as $W(s)$.

   Let us consider the process
  \begin{align}\label{ef3}
Y(t)=&\int_{-\infty}^{t} T(t-s)f(s+\alpha_k,Y(s))ds +\int_{-\infty}^{t} T(t-s)g(s+\alpha_k,Y(s))\mathrm{d} W(s).
\end{align}
Then, analogous to the proof of Lemma \ref{th3.1}, one can see that \eqref{ef3} admits a unique solution $Y^k(t)$ in $\mathbb{U}$.
Note that $Y(t+\alpha_k)$ has the same distribution as $Y^k(t)$ for every $t\in \mathbb{R}$. We claim that $Y(t+\alpha_k)$
converges  to $y_i^k(t)$ with respect to $\|\cdot\|_{B^p}$.
In fact, by \eqref{ef2} and \eqref{ef3} we have
\begin{align}\label{jh3}
 &\|Y^k(\cdot)-\tilde Y(\cdot)\|^p_{B^p} \nonumber\\
 \leq&4^{p-1} \overline{\lim_{l\to \infty}}\frac{1}{2l}\int_{-l}^{l}E \bigg\|\int_{-\infty}^{t} T(t-s)[f(s+\alpha_k,Y^k(s))-\tilde{f}(s,Y^k(s))]\mathrm{d}s\bigg\|^p\dd t \nonumber\\
 &+4^{p-1} \overline{\lim_{l\to \infty}}\frac{1}{2l}\int_{-l}^{l}E\bigg\| \int_{-\infty}^{t} T(t-s)[\tilde{f}(s,Y^k(s))-\tilde{f}(s,Y(s))]\mathrm{d}s \bigg\|^p\dd t\nonumber\\
 &+4^{p-1} \overline{\lim_{l\to \infty}}\frac{1}{2l}\int_{-l}^{l}E\bigg\|\int_{-\infty}^{t} T(t-s)[g(s+\alpha_k,Y^k(s))-\tilde{g}(s,Y^k(s))]\mathrm{d} W(s)\bigg\|^p\dd t\nonumber\\
 &+4^{p-1} \overline{\lim_{l\to \infty}}\frac{1}{2l}\int_{-l}^{l}E\bigg\|\int_{-\infty}^t T(t-s)[\tilde{g}(s,Y^k(s))-\tilde{g}(s,Y(s))]\mathrm{d} W(s)\bigg\|^p\dd t\nonumber\\
 =:&B_1(k,p)+B_2(k,p)+B_3(k,p)+B_4(k,p).
\end{align}

According to the H\"older inequality, the Fubini theorem, the Lebesgue dominated convergence theorem  and \eqref{jh1}, we deduce that
\begin{align}\label{jh4}
&B_1(k,p) \nonumber\\
\leq &4^{p-1}M^p \bigg(\frac{1}{\delta}\bigg)^{p-1} \overline{\lim_{l\to \infty}}\frac{1}{2l}\int_{-l}^{l}\int_{-\infty}^{t} e^{-\delta(t-s)}E\|f(s+\alpha_k,Y^k(s))-\tilde{f}(s,Y^k(s))\|^p\mathrm{d}s \dd t \nonumber\\
\leq &4^{p-1}M^p \bigg(\frac{1}{\delta}\bigg)^{p-1}\int_{0}^{\infty} e^{-\delta s} \overline{\lim_{l\to \infty}}\frac{1}{2l}\int_{-l}^{l}E\|f(t-s+\alpha_k,Y^k(t-s))-\tilde{f}(t-s,Y^k(t-s))\|^p\dd t\dd s \nonumber\\
\leq &4^{p-1}M^p \bigg(\frac{1}{\delta}\bigg)^{p-1}\int_{0}^{\infty} e^{-\delta s} \overline{\lim_{l\to \infty}}\frac{2(l+|s|)}{2l}\times\frac{1}{2(l+|s|)}\nonumber\\
&\times\int_{-l-|s|}^{l+|s|}E\|f(t+\alpha_k,Y^k(t))-\tilde{f}(t,Y^k(t))\|^p\dd t\dd s\nonumber \\
= &4^{p-1}M^p \bigg(\frac{1}{\delta}\bigg)^{p} \|f(\cdot+\alpha_k,Y^k(\cdot))-\tilde{f}(\cdot,Y^k(\cdot))\|_{B^p}^p  \rightarrow0 \text{ as } k\rightarrow\infty,
\end{align}
\begin{align}\label{jh5}
B_2(k,p)\leq & 4^{p-1} M^p \bigg(\frac{1}{\delta}\bigg)^{p-1}\overline{\lim_{l\to \infty}}\frac{1}{2l}\int_{-l}^{l} \int_{-\infty}^{t} e^{-\delta(t-s)}E\|\tilde{f}(s,Y^k(s))-\tilde{f}(s,Y(s))\|^p\mathrm{d}s \dd t\nonumber\\
\leq & 4^{p-1} M^p \bigg(\frac{1}{\delta}\bigg)^{p-1} \int_{0}^{\infty} e^{-\delta s}\overline{\lim_{l\to \infty}}\frac{2(l+|s|)}{2l}\frac{1}{2(l+|s|)}\nonumber\\
&\times \int_{-l-|s|}^{l+|s|}E\|\tilde{f}(t,Y^k(t))-\tilde{f}(t,Y(t))\|^p \dd t\mathrm{d}s\nonumber\\
= & 4^{p-1} M^p \bigg(\frac{1}{\delta}\bigg)^{p} \|\tilde{f}(\cdot,Y^k(\cdot))-\tilde{f}(\cdot,Y(\cdot))\|_{B^p}^p \nonumber\\
\leq & 4^{p-1} M^p L^p_f\bigg(\frac{1}{\delta}\bigg)^{p} \|Y^k(\cdot)-Y(\cdot)\|_{B^p}^p.
\end{align}

For $p> 2$, based on Lemma \ref{bdg}, the H\"older inequality, the Fubini theorem, the Lebesgue dominated convergence theorem, and \eqref{jh1}, we deduce that
\begin{align}
B_3(k,p)
\leq  &4^{p-1}C_p M^p\overline{\lim_{l\to \infty}}\frac{1}{2l}\int_{-l}^{l}E\bigg(\int_{-\infty}^{t} e^{-2\delta(t-s)}\|g(s+\alpha_k,Y^k(s))-\tilde{g}(s,Y^k(s))\|^2_{\mathbb{L}_2^0}\mathrm{d} s\bigg)^{\frac{p}{2}}\dd t\nonumber\\
\leq  &4^{p-1}C_p M^p\bigg(\frac{p-2}{p\delta}\bigg)^{\frac{p-2}{2}}\int_{0}^{\infty} e^{-\frac{p}{2}\delta s}\mathrm{d} s\overline{\lim_{l\to \infty}}\frac{1}{2l}\int_{-l}^{l} E\|g(t-s+\alpha_k,Y^k(t-s))\nonumber\\
&-\tilde{g}(t-s,Y^k(t-s))\|^p_{\mathbb{L}_2^0} \dd t\nonumber\\
\leq  &4^{p-1}C_p M^p\bigg(\frac{p-2}{p\delta}\bigg)^{\frac{p-2}{2}}\int_{0}^{\infty} e^{-\frac{p}{2}\delta s}\mathrm{d} s\overline{\lim_{l\to \infty}}\frac{2(l+|s|)}{2l}\nonumber\\
&\times\frac{1}{2(l+|s|)}\int_{-l-|s|}^{l+|s|} E\|g(t+\alpha_k,Y^k(t))-\tilde{g}(t,Y^k(t))\|^p_{\mathbb{L}_2^0} \dd t \nonumber\\
= &4^{p-1}C_p M^p\bigg(\frac{p-2}{p\delta}\bigg)^{\frac{p-2}{2}}\bigg(\frac{2}{p\delta}\bigg) \|g(\cdot+\alpha_k,Y^k(\cdot)) -\tilde{g}(\cdot,Y^k(\cdot))\|^p_{B^p} \dd t \rightarrow 0 \text{ as } k\rightarrow\infty,\label{jh6}\\
  B_4(k,p)
  \leq &4^{p-1}C_p M^p\overline{\lim_{l\to \infty}}\frac{1}{2l}\int_{-l}^{l}E\bigg(\int_{-\infty}^t e^{-2\delta(t-s)}\|\tilde{g}(s,Y^k(s))-\tilde{g}(s,Y(s))\|^2_{\mathbb{L}_2^0}\mathrm{d} s\bigg)^{\frac{p}{2}}\dd t\nonumber\\
   \leq &4^{p-1}C_p M^p \bigg(\frac{p-2}{p\delta}\bigg)^{\frac{p-2}{2}}\int_{0}^\infty e^{-\frac{p}{2}\delta s}\overline{\lim_{l\to \infty}}\frac{2(l+|s|)}{2l} \times \frac{1}{2(l+|s|)}\nonumber\\
   &\times\int_{-l-|s|}^{l+|s|}E \|\tilde{g}(s,Y^k(s))-\tilde{g}(s,Y(s))\|^p \dd t \mathrm{d} s\nonumber\\
  \leq& 4^{p-1}C_p M^pL_g^p\bigg(\frac{p-2}{p\delta}\bigg)^{\frac{p-2}{2}}\bigg(\frac{2}{p\delta}\bigg)\|Y^k(\cdot)-Y(\cdot)\|_{B^p}^p.\label{jh7}
\end{align}

For $p= 2$, by virtue of the It\^{o} isometry, the H\"older inequality, the Fubini theorem, the Lebesgue dominated convergence theorem, and \eqref{jh1}, we can infer that
\begin{align}
 B_3(k,2)
 \leq &4 M^2 \int_{0}^{\infty} e^{-2\delta s}\overline{\lim_{l\to \infty}}\frac{1}{2l}\int_{-l-|s|}^{l+|s|}E\|g(t+\alpha_k,Y^k(t))-\tilde{g}(t,Y^k(t))\|^2_{\mathbb{L}_2^0}\dd t \mathrm{d}s\nonumber\\
 \leq &4 M^2 \bigg(\frac{1}{2\delta}\bigg)\|g(\cdot+\alpha_k,Y^k(\cdot))-\tilde{g}(\cdot,Y^k(\cdot))\|^2_{B^p}
 \rightarrow 0 \text{ as } k\rightarrow 0,\label{jh8}\\
B_4(k,2)
\leq & 4 M^2\overline{\lim_{l\to \infty}}\frac{1}{2l}\int_{-l}^{l} \int_{-\infty}^t e^{-2\delta(t-s)}E\|\tilde{g}(s,Y^k(s))-\tilde{g}(s,Y(s))\|^2_{\mathbb{L}_2^0}\mathrm{d}s \dd t\nonumber\\
\leq & 4 M^2 \int_{0}^\infty e^{-2\delta s}\overline{\lim_{l\to \infty}}\frac{1}{2l}\int_{-l}^{l} E\|\tilde{g}(t-s,Y^k(t-s))-\tilde{g}(t-s,Y(t-s))\|^2_{\mathbb{L}_2^0} \dd t\mathrm{d}s\nonumber\\
\leq & 4 M^2 \int_{0}^\infty e^{-2\delta s}\mathrm{d}s\overline{\lim_{l\to \infty}}\frac{2(l+|s|)}{2l}\times \frac{1}{2(l+|s|)}\int_{-l-|s|}^{l+|s|} E\|\tilde{g}(s,Y^k(s))-\tilde{g}(s,Y(s))\|^2_{\mathbb{L}_2^0} \dd t\nonumber\\
\leq& 4 M^2L_g^2 \frac{1}{2\delta}\|Y^k(\cdot)-Y(\cdot)\|^2_{B^p}.\label{jh9}
\end{align}

Consequently, for $p>2$, from \eqref{jh3}-\eqref{jh7} it follows that
\begin{align*}
	 \|Y^k(\cdot)-\tilde Y(\cdot)\|^p_{B^p}  \le \Lambda_1(1-R_1)^{-1},
\end{align*}
where $R_1$ is defined in $(H_5)$ and $\Lambda_1:=B_1(k,p)+B_3(k,p)\rightarrow0$ as $k\rightarrow\infty$.  Hence,  we have
\begin{align}\label{cd1}
	\lim_{k\to \infty}\|Y^k(\cdot)-\tilde Y(\cdot)\|_{B^p}^p =0.
\end{align}

Similarly, using \eqref{jh2}, we can get
\begin{align}\label{cd2}
	 \lim_{k\to \infty}\|\tilde Y(\cdot-\alpha_k)- Y^k(\cdot)\|_{B^p}^p   =0.
\end{align}

When $p=2$, by \eqref{jh3}-\eqref{jh5}, \eqref{jh8} and \eqref{jh9}, we have
\begin{align*}
	\lim_{k\to \infty}\|Y^k(\cdot)-\tilde Y(\cdot)\|_{B^p}^2  \le \Lambda_2(1-R_2)^{-1},
\end{align*}
where $R_2$ is given in $(H_5)$ and $\Lambda_2=B_1(k,2)+B_3(k,2)\rightarrow0$ as $k\rightarrow\infty$.
Then, we get
\begin{align}\label{cd3}
	\lim_{k\to \infty} \|Y^k(\cdot)-\tilde Y(\cdot)\|_{B^p}^2   =0.
\end{align}

Similarly, we have
\begin{align}\label{cd4}
\lim_{k\to \infty}\|\tilde Y(\cdot-s_k)- Y^k(\cdot)\|_{B^p}^2 =0.
\end{align}

Since $Y(t+s_k)$ has the same distribution as $Y^k(t)$ for all $t\in \mathbb{R}$, from \eqref{cd1} and \eqref{cd3}  it follows that
\begin{align*}
	&\overline{\lim_{h\to \infty}}\frac{1}{\nu ([-h,h])}\int_{-h}^{h}d_{BL}^p(P \circ [Y(t+\alpha_k)]^{-1}, P \circ [Y(t)]^{-1})\dd \mu (t)\\
=&\overline{\lim_{h\to \infty}}\frac{1}{\nu ([-h,h])}\int_{-h}^{h}d_{BL}^p(P \circ [Y^k(t)]^{-1}, P \circ [Y(t)]^{-1})\dd \mu (t)\\
	\le &\overline{\lim_{h\to \infty}}\frac{1}{\nu ([-h,h])}\int_{-h}^{h} \bigg(\int_\Omega  \|Y^k(t)-Y(t)\|\dd P \bigg)^p\dd \mu (t)\nonumber\\
	\leq&\overline{\lim_{h\to \infty}}\frac{1}{\nu ([-h,h])}\int_{-h}^{h}  E\|Y^k(t)-Y(t)\|^p \dd \mu (t)\\
=&\|Y^k(\cdot)-\tilde Y(\cdot)\|_{B^p}^p \rightarrow 0 \text{ as } k\rightarrow \infty.
\end{align*}

Similarly, in view of \eqref{cd2} and \eqref{cd4}, we can get
\begin{align*}
	\lim_{k\to \infty}\overline{\lim_{h\to \infty}}\frac{1}{\nu ([-h,h])}\int_{-h}^{h}d_{BL}^p(P \circ [\tilde{Y}(t-\alpha_k)]^{-1}, P \circ [Y(t)]^{-1})\dd \mu (t)=0.
\end{align*}

  Thus,  $Y$ is $p$-th Besicovitch almost automorphic in distribution. The proof is finished.
\end{proof}

Based on Lemma $\ref{th3.1}$ and Theorem $\ref{th3.2}$,  we will prove the following result.
\begin{theorem}\label{thm3}
Let  $\mu,\nu\in\mathcal{M}_{p}$ and conditions $(A_1)$ and $(A_2)$ hold.	If $(H_1)$, $(H_3)-(H_6)$ are satisfied, then system \eqref{1} has a unique Besicovitch-like $(\mu,\nu )$-pseudo almost automorphic solution in distribution.
\end{theorem}
\begin{proof}
		According to Lemma $\ref{th3.1}$,  system $\eqref{1}$ admits a unique   solution $X\in \mathbb{U}$ satisfying \eqref{eq1}. Moreover, by Theorem $\ref{th3.2}$, we see that  the following system
\begin{align}\label{ef1}
Y(t)=\int_{-\infty}^t T(t-s)f_1(s,Y(s))ds+\int_{-\infty}^t T(t-s)g_1(s,Y(s))\mathrm{d} W(s)
\end{align}
 has a unique Besicovitch almost automorphic solution $Y\in \mathbb{U}$ in distribution.
		
Now, let $X-Y=Z$, then $X=Y+Z$. To complete the proof, we just need to prove that $Z \in\varepsilon_ D(\mathbb{R},\mathcal{L}^p(\Omega,\mathbb{X}),\mu,\nu)$. From \eqref{eq1} and \eqref{ef1}, it follows that
\begin{align}\label{w1}
  E\|Z(t)\|^p
  \leq &4^{p-1}E\bigg\|\int_{-\infty}^t T(t-s)[f(s,X(s))-f_1(s,X(s))]\mathrm{d}s\bigg\|^p\nonumber\\
  &+4^{p-1}E\bigg\|\int_{-\infty}^t T(t-s)[f_1(s,X(s))-f_1(s,Y(s))]\mathrm{d}s\bigg\|^p\nonumber\\
  &+4^{p-1}E\bigg\|\int_{-\infty}^t T(t-s)[g(s,X(s))-g_1(s,X(s))]\mathrm{d} W(s)\bigg\|^p\nonumber\\
  &+4^{p-1}E\bigg\|\int_{-\infty}^t T(t-s)[g_1(s,X(s))-g_1(s,Y(s))]\mathrm{d} W(s)\bigg\|^p\nonumber\\
  =:&D_1(t,p)+D_2(t,p)+D_3(t,p)+D_4(t,p).
\end{align}
	
Then, we have
\begin{align}
 D_1(t,p)
 \leq& 4^{p-1}M^p\bigg(\frac{1}{\delta}\bigg)^{p-1}\int_{-\infty}^t e^{-\delta(t-s)}E\|f_2(s,X(s))\|^p \mathrm{d}s=:\varrho_1(t,p),\label{w2}\\
 D_2(t,p)
 \leq &4^{p-1}M^pL_f^p\bigg(\frac{1}{\delta}\bigg)^{p-1}\int_{-\infty}^t e^{-\delta(t-s)}E\|Z(s)\|^p\mathrm{d}s=:\varrho_2(t,p).\label{w3}
\end{align}

When $p>2$, by Lemma \ref{bdg} and the Fubini theorem, we can infer that
\begin{align}
 D_3(t,p)
 \leq & 4^{p-1}C_p M^p \bigg(\frac{p-2}{p\delta}\bigg)^{\frac{p-2}{2}}\int_{-\infty}^t e^{-\frac{p}{2}\delta(t-s)}E\|g_2(s,X(s))\|^p\mathrm{d}s=:\varrho_3(t,p),\label{w4}\\
 D_4(t,p)
 \leq &4^{p-1}C_p M^p L_g^p\bigg(\frac{p-2}{p\delta}\bigg)^{\frac{p-2}{2}}\int_{-\infty}^t e^{-\frac{p}{2}\delta(t-s)}E\|Z(s)\|^p\mathrm{d}s=:\varrho_4(t,p),\label{w5}
\end{align}
and, when $p=2$, according to the It\^{o} isometry and the Fubini theorem, we obtain
\begin{align}
 D_3(t,2)
 \leq & 4   M^2  \int_{-\infty}^t e^{-2\delta(t-s)}E\|g_2(s,X(s))\|^2\mathrm{d}s=:\hat{\varrho}_3(t,2),\label{w6}\\
 D_4(t,2)
 \leq &4 M^2 L_g^2 \int_{-\infty}^t e^{-2\delta(t-s)}E\|Z(s)\|^2\mathrm{d}s=:\hat{\varrho}_4(t,2).\label{w7}
\end{align}

From  \eqref{w1}-\eqref{w5}, it follows that for $p>2$,
	\begin{align*}
		 E\|Z(t)\|^p
		\le \chi_1(t,p)+\beta_1\int_{-\infty}^{t}e^{-\delta (t-s)} E\|Z(s)\|^p \dd s,
	\end{align*}
	where $\chi_1(t,p):=\varrho_1(t,p)+\varrho_3(t,p)$, and $\beta_1$ is defined in $(H_6)$.

In view of Lemma \ref{GR}, for any $\gamma\in (0,\delta-\beta_1)$, there holds
\begin{align}\label{yy1}
E\|Z(t)\|^p\leq \chi_1(t,p)+\beta_1\int_{-\infty}^{t}e^{-\gamma(t-s)}\chi_1(s,p)\dd s.
\end{align}

Note that
\begin{align*}
&  \lim_{h\to \infty}\frac{1}{\nu ([-h,h])}\int_{-h}^{h}\int_{-\infty}^t e^{-\delta(t-s)}E\|f_2(s,X(s))\|^p \mathrm{d}s \dd \mu (t)\\
=&\lim_{h\to \infty}\frac{1}{\nu ([-h,h])}\int_{-h}^{h}\int_{0}^{\infty}e^{-\delta s}  E\|f_2(t-s,X(t-s)\|^p\dd s \dd \mu (t)\\
=&\int_{0}^{\infty}e^{-\delta s}\lim_{h\to \infty}\frac{1}{\nu ([-h,h])}\int_{-h}^{h}   E\|f_2(t-s,X(t-s)\|^p \dd \mu (t)\dd s\\
\leq &\int_{0}^{\infty}e^{-\delta s} \frac{\nu_s ([-h-s,h+s])}{\nu ([-h,h])}\times \frac{1}{\nu_s ([-h-s,h+s])}\int_{-h-s}^{h+s}   E\|f_2(t,X(t)\|^p \dd \mu_s (t)\dd s=0
\end{align*}
and
\begin{align*}
&  \lim_{h\to \infty}\frac{1}{\nu ([-h,h])}\int_{-h}^{h}\int_{-\infty}^{t}e^{-\gamma (t-s)}\int_{-\infty}^{s}e^{-\delta(s-u)} E\|f_2(u,X(u))\|^p\dd u \dd s\dd \mu (t)\\
=& \lim_{h\to \infty}\frac{1}{\nu ([-h,h])}\int_{-h}^{h}\int_{0}^{\infty}e^{-\gamma  s}\int_{-\infty}^{t-s}e^{-\delta(t-s-u)} E\|f_2(u,X(u))\|^p\dd u\dd s \dd \mu (t)\\
=& \lim_{h\to \infty}\frac{1}{\nu ([-h,h])}\int_{-h}^{h}\int_{0}^{\infty}e^{-\gamma  s}\int_{0}^{\infty}e^{-\delta u} E\|f_2(t-s-u,X(t-s-u))\|^p\dd u\dd s \dd \mu (t)\\
=&\int_{0}^{\infty}e^{-\gamma  s}\int_{0}^{\infty}e^{-\delta u}\lim_{h\to \infty}\frac{1}{\nu ([-h,h])}\int_{-h}^{h} E\|f_2(t-s-u,X(t-s-u))\|^p\dd \mu (t)\dd u\dd s\\
\leq&\int_{0}^{\infty}e^{-\gamma s}\int_{0}^{\infty}e^{-\delta u}\lim_{h\to \infty}\frac{1}{\nu ([-h,h])}\int_{-h-s-u}^{h+s+u}  E\|f_2(t,X(t))\|^p\dd \mu_{s+u} (t)\dd u\dd s=0.
\end{align*}
We can get
\begin{align*}
\lim_{h\to \infty}\frac{1}{\nu ([-h,h])}\int_{-h}^{h}\varrho_1(t,p)\dd \mu (t)=\lim_{h\to \infty}\frac{1}{\nu ([-h,h])}\int_{-h}^{h}\int_{-\infty}^{t}e^{-\gamma(t-s)}\varrho_1(s,p)\dd s\dd \mu(t)=0.
\end{align*}

Similarly,   we can obtain
\begin{align*}
\lim_{h\to \infty}\frac{1}{\nu ([-h,h])}\int_{-h}^{h}\varrho_3(t,p)\dd \mu (t)=\lim_{h\to \infty}\frac{1}{\nu ([-h,h])}\int_{-h}^{h}\int_{-\infty}^{t}e^{-\gamma(t-s)}\varrho_3(s,p)\dd s\dd \mu(t)=0.
\end{align*}

Consequently, by \eqref{yy1}, we arrive at
\begin{align*}
\lim_{h\to \infty}\frac{1}{\nu ([-h,h])}\int_{-h}^{h}E\|Z(t)\|^p\dd \mu (t)=0.
\end{align*}

	When $p=2$,  based on \eqref{w1}-\eqref{w3}, \eqref{w6}, \eqref{w7}, and Lemma \ref{GR}, for any $\check{\gamma}\in (0,\delta-\beta_2)$, one can obtain that
	\begin{align}\label{yy3}
		 E\|z(t)\|_0^2
		\le  \chi_2(t,2)(t)+\beta_2 \int_{-\infty}^{t}e^{-\check{\gamma} (t-s)} \chi_2(s,2) \dd s,
	\end{align}
	where $\chi_2(t,2):=\varrho_1(t,2)+\hat{\varrho}_3(t,2)$ and $\beta_2$ is mentioned in $(H_6)$.

Similar to the proof of the case for $p>2$, from \eqref{yy3}, one can get
	\begin{align*}
		\lim_{h\to \infty}\frac{1}{\nu ([-h,h])}\int_{-h}^{h}E\|Z(t)\|^2\dd \mu (t)=0.
	\end{align*}

Hence,  $Z \in\varepsilon_ M(\mathbb{R},\mathcal{L}^p(\Omega,\mathbb{X}),\mu,\nu)$. Consequently, we obtain that   $Z \in\varepsilon_ D(\mathbb{R},\mathcal{L}^p(\Omega,\mathbb{X}),\mu,\nu)$ by Theorem \ref{key}. The proof is finished.
\end{proof}

\section{Global exponential stability}
\setcounter{equation}{0}
\indent

This section is devoted to the global exponential stability of solutions to system \eqref{1}.

\begin{theorem}\label{thm2}
Under the assumptions of Theorem \ref{thm3}, let $X$ be the Besicovitch-like $(\mu,\nu)$-pseudo almost automorphic solution in distribution and $Y$ be any solution. Then there exist $\lambda>0$ and $N>0$ such that
\begin{align}\label{st}
  E\|Y(t)-X(t)\|^{p} \le N\,E\|Y(t_0)-X(t_0)\|^{p} \,\, e^{-p\lambda (t-t_0)},\quad t>t_0.
\end{align}
\end{theorem}

\begin{proof} Let $Y(t)-X(t)=Z(t)$, then from \eqref{jf1}, we have
\begin{align}\label{ss}
 Z(t)=&T(t-t_0)Z(t_0)+\int_{t_0}^tT(t-s)[f(s,Y(s))-f(s,X(s))\dd s\nonumber\\
 &+\int_{t_0}^t T(t-s)[g(s,Y(s))-g(s,X(s))]\dd W(s).
\end{align}

For $p>2$, let us consider the  linear function $\Theta: \mathbb{R}\rightarrow \mathbb{R}$ defined as follows:
\begin{align*}
 \Theta(\theta)= \delta-p\theta- 3^{p-1}M^p \Xi,
\end{align*}
where $\Xi:=L_f^p\big(\frac{1}{\delta}\big)^{p-1}+ C_p   L_g^p\big(\frac{p-2}{p\delta}\big)^{\frac{p-2}{2}}$.

Obviously, by condition $(H_6)$, there exists a unique $\theta^*>0$ such that $\Theta(\theta^*)=0$ and $\Theta(\theta)>0$ for all $\theta\in (0,\theta^*)$. Choose a $\lambda \in (0, \min\{\theta^*,\delta/p\})$, then $\Theta(\lambda)>0$,  i.e.,
 \begin{align*}
 3^{p-1}M^p\Xi(\delta-p\lambda)^{-1}<1,
 \end{align*}

Let $N= \frac{\delta}{\Xi}$, then
\begin{align*}
 \frac{1}{N}-\frac{\Xi}{\delta-p\lambda}<0.
\end{align*}

We claim that for any $\varepsilon>0$, it holds
		\begin{align}\label{s1}
			\mathbb{E}\|Z(t)\|^p\leq N (\mathbb{E}\|Z(t_0)\|^p + \varepsilon) e^{-p\lambda(t-t_0)}, \quad t \geq t_0.
		\end{align}
		
Otherwise, there exists a $\bar{t} > 0$ such that
		\begin{equation}\label{s2}\left\{\begin{array}{ll}
			E\|Z(\bar{t})\|^p= NE\|Z(t_0)\|^p+ \varepsilon) e^{-p\lambda(\bar{t}-t_0)},\\
			E \|Z(t)\|^p \leq N E\|Z(t_0)\|^p + \varepsilon ) e^{-p\lambda(t-t_0)}, \quad t_0\leq t<\bar{t}.
\end{array}\right.
		\end{equation}
		
Based on \eqref{ss} and \eqref{s2}, we infer that
\begin{align*}
&E \|Z(\bar{t})\|^p\\
 \leq&3^{p-1}M^p e^{-p\delta(\bar{t}-t_0)}E\|Z(t_0)\|^p+3^{p-1}M^pL_f^p\bigg(\frac{1}{\delta}\bigg)^{p-1}\int_{t_0}^{\bar{t}}e^{-\delta(\bar{t}-s)}E\|Z(s)\|^p\dd s\\
 &+3^{p-1}C_p M^p L_g^p \bigg(\frac{p-2}{p\delta}\bigg)^{\frac{p-2}{2}} \int_{t_0}^{\bar{t}} e^{-\frac{p}{2}\delta(\bar{t}-s)}E\|Z(s)\|^p\dd s\\
 \leq&3^{p-1}M^p e^{-p\delta(\bar{t}-t_0)}(E\|Z(t_0)\|^p+\varepsilon)+3^{p-1}M^pL_f^p\bigg(\frac{1}{\delta}\bigg)^{p-1}N (E\|Z(t_0)\|^p + \varepsilon)\int_{t_0}^{\bar{t}}e^{-\delta(\bar{t}-s)}\\
 &\times e^{-p\lambda(s-t_0)}\dd s+3^{p-1}C_p M^p L_g^p \bigg(\frac{p-2}{p\delta}\bigg)^{\frac{p-2}{2}}N (E\|Z(t_0)\|^p + \varepsilon)\int_{t_0}^{\bar{t}} e^{-\frac{p}{2}\delta(\bar{t}-s)} e^{-p\lambda(s-t_0)}\dd s\\
 \leq&N(E\|Z(t_0)\|^p+\varepsilon) e^{-p\lambda(\bar{t}-t_0)}\bigg[ \frac{3^{p-1}M^p e^{-p(\delta-\lambda)(\bar{t}-t_0)}}{N}
 +3^{p-1}M^pL_f^p\bigg(\frac{1}{\delta}\bigg)^{p-1} \\
 &\times\frac{1}{\delta-p\lambda} \big(1-e^{-(\delta-p\lambda)(\bar{t}-t_0)}\big)+3^{p-1}C_p M^p L_g^p \bigg(\frac{p-2}{p\delta}\bigg)^{\frac{p-2}{2}}  \frac{1}{\frac{p}{2}\delta-p\lambda} \big(1-e^{-(\frac{p\delta}{2}-p\lambda)(\bar{t}-t_0)}\big) \bigg]\\
 \leq&N(E\|Z(t_0)\|^p+\varepsilon) e^{-p\lambda(\bar{t}-t_0)}\bigg[ \frac{3^{p-1}M^p e^{-p(\delta-\lambda)(\bar{t}-t_0)}}{N}
 +3^{p-1}M^pL_f^p\bigg(\frac{1}{\delta}\bigg)^{p-1}\frac{1}{\delta-p\lambda} \\
 &\times\big(1-e^{-p(\delta-\lambda)(\bar{t}-t_0)}\big)+3^{p-1}C_p M^p L_g^p \bigg(\frac{p-2}{p\delta}\bigg)^{\frac{p-2}{2}}  \frac{1}{\delta-p\lambda} \big(1-e^{-p(\delta-\lambda)(\bar{t}-t_0)}\big)\bigg] \\
 \leq&N(E\|Z(t_0)\|^p+\varepsilon) e^{-p\lambda(\bar{t}-t_0)}\bigg\{3^{p-1}M^p e^{-p(\delta-\lambda)(\bar{t}-t_0)}\bigg[\frac{1}{N}-\frac{\Xi}{\delta-p\lambda}\bigg]+\frac{3^{p-1}M^p\Xi}{\delta-p\lambda}
  \bigg\}\\
  <&N(E\|Z(t_0)\|^p+\varepsilon) e^{-p\lambda(\bar{t}-t_0)},
\end{align*}
which contradicts the first formula in \eqref{s2}. Hence,   \eqref{s1} holds.
		 Letting $\varepsilon \to 0^{+}$, we obtain that \eqref{st} holds for $p>2$.

For $p=2$,  consider the  linear function $\Theta: \mathbb{R}\rightarrow \mathbb{R}$ defined as follows:
\begin{align*}
 \tilde{\Theta}(\theta)= \delta-2\theta- 3 M^2 \hat{\Xi},
\end{align*}
where $\Xi:=L_f^p\big(\frac{1}{\delta}\big) + L_g^2$. Similarly, one can show that \eqref{st} also holds for $p=2$.

   This completes the proof.
	\end{proof}

\section{Example}
\indent

\setcounter{equation}{0}
 \indent

Consider the following stochastic partial differential equation:
\begin{equation}\label{ex1}
\left\{\begin{aligned}
&\mathrm{d} u(t,\xi)=\displaystyle \frac{\partial^2}{\partial \xi^2} u(t,\xi)+F(t,u(t,\xi))\mathrm{d}t+G(t,u(t,\xi) )\mathrm{d}W(t), \,\, \text{ for } (t,\xi) \in (0,\pi)\times \mathbb{R},\\
&u(t,0)=u(t,\pi)=0, \,\, \mathrm{for}\,\,  t\in \mathbb{R},
\end{aligned}\right.
\end{equation}
where $W(t)$ is a one-dimensional two-sided Brownian motion,
\begin{align*}
 & F(t,u(t,\xi))\\
  =&\frac{1}{60}  \arctan t+\frac{1}{30}\left\{\frac{1}{3}\left[\cos\left(\frac{1}{2+\sin t+\sin\sqrt{5}t}\right)+e^{-|t|}\right] -\frac{1}{2}\arctan t\right\} \tanh u(t,\xi)\\
&G(t,u(t,\xi)),\\
=&\frac{1}{60}\bigg[\cot^{-1} t+\frac{1}{2}\left(\sin t- \cos \sqrt{3}t
 +\frac{\cos t}{1+t^2}+\mathrm{sign}(t)\right) \tanh u (t,\xi)+\frac{1}{2}\cot^{-1} t\cos u (t,\xi)\bigg],
\end{align*}

Take $\mathbb{X} = \mathcal{L}^2[0, \pi]$ with inner product $\langle\cdot,\cdot\rangle$ and norm $\|\cdot\|$. The operator $A : D(A)\subset \mathbb{X}\rightarrow\mathbb{X}$ is the Laplace
operator, where
$D(A)=H^2(0,\pi)\cap H^1_0(0,\pi).$

It is well known that  operator $A$ generates  a $C_0$-semigroup $T(t)$ on $\mathbb{X}$ with $\|T(t)\|\leq e^{- t}$ for $t\geq 0$. Hence, condition $(H_1)$ satisfied with $M=1$ and $\delta=1$.

Then system \eqref{ex1} can be rewritten as an abstract  equation
\begin{align*}
\mathrm{d}x(t)=&\big[Ax(t)+f(t,x(t)))\big]\mathrm{d}t+g(t,x(t))\mathrm{d}W(t),\,\, t\in\mathbb{R},
\end{align*}
on the Hilbert space $\mathbb{X}$, where $x(t):=u(t,\cdot),$
\begin{align*}
  f(t,x(t))=F(t,u(t,\xi)) \text{ and }
 g(t,x(t))=G(t,u(t,\xi))).
\end{align*}

Let
\begin{align*}
f(t,x(t))=f_1(t,x(t))))+f_2(t,x(t)) \text{ and }
 g(t,x(t))=g_1(t,x(t))+g_2(t,x(t)),
\end{align*}
where
\begin{align*}
  f_1(t,x(t))=&\frac{1}{90}\left[\cos\left(\frac{1}{2+\sin t+\sin\sqrt{5}t}\right)\right] \tanh x(t),\\
  f_2(t,x(t))=&\frac{1}{60}\arctan t-\frac{1}{60} \arctan t\tanh x(t) +\frac{1}{90}e^{-|t|}\tanh x(t),\\
  g_1(t,x(t))=&\frac{1}{120}\left(\sin t- \cos \sqrt{3}t
 \right)\tanh x (t),\\
  g_2(t,x(t))=&\frac{1}{120}\left[(2+\cos x (t))\cot^{-1} t+\bigg(\frac{\cos t}{1+t^2} +\mathrm{sign}(t)\bigg)\tanh x(t)\right].
\end{align*}

Take the Radon-Nikodym derivatives of $\mu$ and $\nu$ are $2\cos^2 t+3$ and $e^{|t|}+\cos t+2$ with respect to the Lebesgue measure on $\mathbb{R}$, respectively. Then, we have
\begin{align*}
 \overline{\lim_{h\to \infty}} \frac{\mu([-h,h])}{\nu([-h,h])}=\overline{\lim_{h\to \infty}} \frac{\int_{-h}^{h}(2\cos^2t+3)\mathrm{d}t}{\int_{-h}^{h}(e^{|t|}+\cos t+2)\mathrm{d}t}<\infty,
\end{align*}
and hence, condition $(A_1)$ is satisfied. Moreover, for all $\sigma\in \mathbb{R}$ and $A\in \Sigma$, we have
\begin{align*}
  \mu_\sigma(A)=\int_{A}[2\cos^2(t-\sigma)+3]\mathrm{d}t\leq \frac{5}{3}\int_{A}[2\cos^2t+3]\mathrm{d}t\leq \frac{5}{3} \mu (A).
\end{align*}
Take $\rho=\frac{5}{3}$ and $I=\emptyset$,
then  condition  $(A_2)$ is also satisfied.

 Since   $\cos\big(\frac{1}{2+\sin t+\sin\sqrt{5}t}\big)$ is almost automorphic in the Bochner sense  and $\sin t- \cos \sqrt{3}t$ is almost periodic in the Bohr sense,
   $f_1\in B^pAA_m(\mathbb{R}\times (\mathcal{L}^{p}(\Omega,\mathbb{X})), \mathcal{L}^{p}(\Omega,\mathbb{X}))$ and $g_1\in B^pAA_m(\mathbb{R}\times (\mathcal{L}^{p}(\Omega,\mathbb{H}))^2, \mathcal{L}^{p}(\Omega,\mathbb{L}_{2}^{0}))$.

Note that
\begin{align*}
&  \lim\limits_{h\rightarrow +\infty}\frac{1}{\nu([-h,h])}\int_{-h}^{h} E\left\| f_2(t,x(t)))\right\|^p\mathrm{d}\mu(t)
  \leq\left(\frac{\pi}{60}+\frac{1}{90}\right)\lim\limits_{h\rightarrow\infty}\frac{\int_{-h}^{h}(2\cos^2t+3)\mathrm{d}t}{\int_{-h}^{h}(e^{|t|}+\cos t+2)\mathrm{d}t}=0,\\
&  \lim\limits_{h\rightarrow +\infty}\frac{1}{\nu([-h,h])}\int_{-h}^{h} E\left\| g_2(t,x(t))\right\|_{\mathbb{L}_2^0}^p\mathrm{d}\mu(t)
  \leq\frac{3\pi+2}{120}\lim\limits_{h\rightarrow\infty}\frac{\int_{-h}^{h}(2\cos^2t+3)\mathrm{d}t}{\int_{-h}^{h}(e^{|t|}+\cos t+2)\mathrm{d}t}=0.
\end{align*}

Thus,
$f_2\in \varepsilon_M(\mathbb{R}\times \mathcal{L}^p(\Omega,\mathbb{X}),\mathcal{L}^p(\Omega,\mathbb{X}),\mu,\nu)$ and $g_2\in \varepsilon_M(\mathbb{R}\times \mathcal{L}^p(\Omega,\mathbb{H}),\mathcal{L}^{p}(\Omega,\mathbb{L}_{2}^{0}),\mu,\nu)$.

Take  $L_f=\frac{1}{90}+\frac{\pi}{120}, L_g=\frac{4+\pi}{120}, M_f=\big(\frac{1}{45}+\frac{\pi}{60}\big)^{p}, M_g= \big(\frac{1}{30}+\frac{\pi}{40}\big)^{p}$, then by Remark \ref{re11}, condition $(H_4)$ is satisfied.

By computing, we get
  \begin{align*}
  \beta_1:=&4^{p-1}M^p\bigg[L_f^p\bigg(\frac{1}{\delta}\bigg)^{p-1}+ C_p   L_g^p\bigg(\frac{p-2}{p\delta}\bigg)^{\frac{p-2}{2}}\bigg]\approx 0.09091<1=\delta, \,\, (p=\frac{7}{2}),\\
  \beta_2:=&4 M^2\bigg[L_f^2\bigg(\frac{1}{\delta}\bigg) + L_g^2\bigg]\approx 0.01973<1=\delta.
\end{align*}

Hence, for $p=\frac{7}{2}$ and $p=2$, condition $(H_5)$ are satisfied.

Consequently, according to Theorem \ref{thm3}, system \eqref{ex1} admits a unique $(\mu,\nu)$ pseudo $p$-th Besicovitch-like almost automorphic mild solution in distribution that is globally exponentially stable.
	
\begin{remark}
In the above example, we can also decompose $f$ and $g$ as $f=\check{f}_1+\check{f}_2$ and $g=\check{g}_1+\check{g}_2$, respectively, where
\begin{align*}
  \check{f}_1(t,x(t))=&\frac{1}{90}\left[\cos\left(\frac{1}{2+\sin t+\sin\sqrt{5}t}\right)+e^{-|t|}\right] \tanh x(t),\\
  \check{f}_2(t,x(t))=&\frac{1}{60}\left(1-\tanh x(t)\right)\arctan t,\\
  \check{g}_1(t,x(t))=&\frac{1}{120}\left(\sin t- \cos \sqrt{3}t
 +\frac{\cos t}{1+t^2}+\mathrm{sign}(t)\right)\tanh x (t),\\
  \check{g}_2(t,x(t))=&\frac{1}{120}[2+\cos x (t)]\cot^{-1} t.
\end{align*}

Since  $\|e^{-|t|}\|_{B^p}=0$, $\|\frac{\cos t}{1+t^2}\|_{B^p}=0$, and $\|\mathrm{sign}(\cdot)\|_{B^p}=0$. Hence,  $f_1\in B^pAA_m(\mathbb{R}\times (\mathcal{L}^{p}(\Omega,\mathbb{X})),
 \mathcal{L}^{p}(\Omega,\mathbb{X}))$ and $g_1\in B^pAA_m(\mathbb{R}\times (\mathcal{L}^{p}(\Omega,\mathbb{H}))^2, \mathcal{L}^{p}(\Omega,\mathbb{L}_{2}^{0}))$.
In addition, it is easy to see that $\check{f}_1\in B^pAA_m(\mathbb{R}\times (\mathcal{L}^{p}(\Omega,\mathbb{X})), \mathcal{L}^{p}(\Omega,\mathbb{X}))$ and $\check{f}_2\in \varepsilon_M(\mathbb{R}\times \mathcal{L}^p(\Omega,\mathbb{X}),\mathcal{L}^p(\Omega,\mathbb{X}),\mu,\nu)$.  Therefore, the decomposition of $f$ into an almost automorphic part and an ergodic part is not unique.
\end{remark}

\section{Conclusion}
\indent	

In this paper, we introduced the concept of Besicovitch-like  pseudo almost automorphic stochastic processes. As an application, we established the existence, uniqueness, and global exponential stability of Besicovitch-like pseudo almost automorphic solutions for a class of semilinear stochastic differential equations. The results and methods developed in this paper can also be applied to study the existence of such solutions for other types of stochastic differential equations, including fractional stochastic differential equations, stochastic differential equations driven by fractional Brownian motion, and stochastic differential equations with L\'{e}vy jumps.

\section*{Data availability}
\indent

 No datasets were generated or analysed during the current study.

\section*{Conflicts of Interest}
\indent

The author declares that there are no conflicts of interest.

	{}

\end{document}